\documentclass[12pt]{article}
\usepackage{epsfig}
\usepackage{pdflscape}
\usepackage{amsmath,amssymb}
\usepackage{graphicx}
\usepackage{color}
\usepackage{url}
\usepackage{longtable}
\usepackage[width=\textwidth]{caption}

\definecolor{blue}{rgb}{0,0,0.9}
\definecolor{red}{rgb}{0.9,0,0}
\definecolor{green}{rgb}{0,0.9,0}
\definecolor{purple}{rgb}{0.45,0,0.45}

\newenvironment{Proof}{\noindent {\it Proof. }}{\hfill $\square$}

\newlength{\len}
\newtheorem{theorem}{Theorem}[section]

\newtheorem{remark}{Remark}[section]
\newtheorem{defn}{Definition}[section]

\newcommand{\mb}[1]{\mbox{\boldmath $#1$}}
  \def\h{\mb{h}}

\def\bD{\mb{D}}
\def\bU{\mb{U}}
\def\bQ{\mb{{\cal Q}}}

\def\bB{\mb{{\cal B}}}
 \def\bX{\mb{{\cal X}}} \def\bY{\mb{{\cal Y}}}
\def\bK{\mb{{\cal K}}} \def\bW{\mb{{\cal W}}}
\def\bD{\mbox{\boldmath{$D$}}}

\def\cA{{\cal A}} \def\cB{{\cal B}}  \def\cD{{\cal D}}
\def\cE{{\cal E}}  \def\cF{{\cal F}}  \def\cG{{\cal G}}  
\def\cJ{{\cal J}}  \def\cK{{\cal K}} \def\cL{{\cal L}}  \def\cN{{\cal N}}
 \def\cQ{{\cal Q}}
  \def\cT{{\cal T}}    \def\cV{{\cal V}} \def\cW{{\cal W}}
\def\cX{{\cal X}} \def\cY{{\cal Y}}

\def\alp{\alpha} 
 \def\sig{\sigma}

\def\approxleq{ \kern3pt \mbox{\raisebox{.6ex}{$<$}} \kern-8pt
\mbox{\raisebox{-.6ex}{$\sim$}} \kern5pt}

\def\norm#1{\|#1 \|}
\def\abs#1{|#1 |}
\def\inprod#1#2{\langle#1,\,#2 \rangle}

\def\diag#1{{\rm diag}(#1)}

\def\skron{{{\kern2.3pt\bigcirc\kern-8.5pt\ast\kern5pt }}}

\def\nn{\nonumber}

\begin{document}

\title{\bf A semi-proximal augmented Lagrangian based decomposition method for primal
 block angular convex composite quadratic conic programming problems}

\author{
Xin-Yee Lam
\thanks{Department of Mathematics, National University of
         Singapore, 10 Lower Kent Ridge Road, Singapore 119076
         ({\tt mattohkc@math.nus.edu.sg}).},
\quad
Defeng Sun
\thanks{Department of Applied Mathematics, The Hong Kong Polytechnic University, Hung Hom, Hong Kong ({\ttfamily defeng.sun@polyu.edu.hk}).},
\quad
 Kim-Chuan Toh
 \thanks{Department of Mathematics, and Institute of Operations Research and Analytics, National University of
         Singapore, 10 Lower Kent Ridge Road, Singapore 119076
         ({\tt mattohkc@math.nus.edu.sg}).
         This author's research is supported in part by
the Ministry of Education, Singapore, Academic Research Fund under Grant R-146-000-257-112.
      }
}

%\date{Last updated: December 4, 2018}

\maketitle

\begin{abstract}
We propose a semi-proximal augmented Lagrangian based decomposition method for convex composite quadratic conic programming problems with  primal
 block angular structures.
 %Our method is particularly well suited for solving quadratic conic programming problems arising from stochastic programming models.
Using our algorithmic framework, we are able to naturally
derive several well known augmented Lagrangian based decomposition methods
for stochastic programming such as the diagonal quadratic approximation method of
Mulvey and Ruszczy\'{n}ski. Moreover, we are able to derive novel enhancements and generalizations
of these well known methods.
We also propose a semi-proximal symmetric Gauss-Seidel based alternating direction method of multipliers for solving the corresponding dual problem.
Numerical results show that our algorithms can perform well even for very large instances
of primal block angular convex QP problems. For example, one instance
with more than $300,000$  linear constraints and $12,500,000$ nonnegative variables is solved in less
than a minute whereas Gurobi took more  than 3 hours, and
another instance {\tt qp-gridgen1}  with more than $331,000$ linear constraints
and $986,000$ nonnegative variables is solved in about 5 minutes whereas Gurobi took more than 35 minutes.

\end{abstract}

%%%%%%%%%%%%%%%%%%%%%%%%%%%%%%%%%%%%%
\section{Introduction}

In this paper, we will focus on solving convex composite quadratic conic programming problems with a primal block angular structure, i.e. optimization problems with a separable convex objective function and conic constraints but the variables are coupled by linking linear constraints across different variables. Without specially designed  strategies to exploit
the underlying block angular structure, computational inefficiency of an algorithm will be severe because the constraints cannot be decomposed completely.

In practical applications, quadratic and linear problems with primal block angular structure appear in many contexts, such as multicommodity flow problems \cite{MCPsurvey} and statistical disclosure control \cite{SDCbook}. These problems are often very large scale in practice, and
 standard interior point methods such as those implemented in
 Gurobi or Mosek may not be efficient enough to solve such problems.
In the literature, specialized algorithms designed to solve these problems have been studied extensively. Three of the most widely known algorithmic classes are (i) decomposition methods based on augmented Lagrangian and proximal-point algorithms, see for example
\cite{MR92,RocWet91,Rusz86,Rusz89,ruszcz95,ruszcz99};
 (ii) interior-point log-barrier Lagrangian decomposition methods such as those studied in \cite{zhaogy991,zhao2,zhaogy992,MO07,MO09}; and (iii)
standard interior-point methods which incorporate novel numerical linear algebraic techniques to
exploit the underlying block angular structures
when solving the large linear systems arising in each iteration, for example in
\cite{BirQi88,CG93,GSV97,SM91,Todd88}.

Besides quadratic and linear problems, semidefinite programming (SDP) problems with primal block angular structures are beginning to appear in the literature more frequently. It is gaining more attention as practitioners become more sophisticated in using SDP to model their application problems. For example, the authors in \cite{HK2017} reformulated a two-stage distributionally robust linear program as a completely positive cone program which bears a block angular structure and applied the reformulation to solve a multi-item newsvendor problem. Although linear programming problems with primal block angular structures have been studied extensively, the more complicated SDP version is still in its infancy stage. Apart from \cite{MO07}, \cite{Siv10} and \cite{ZA11}, we are not aware of other works.
%In fact, linear programming could be viewed as a special case of SDP with diagonal matrices.

By focusing on designing efficient algorithms for solving general conic programming problems
with primal block angular structures, we can in general
also use the same algorithmic framework to
solve the primal block angular linear and quadratic programming problems efficiently
through designing novel numerical linear algebraic techniques to exploit the underlying
structures.
In this paper, our main objective is to design efficient and robust (distributed) algorithms for solving large scale conic programming problems with block angular structures. Specifically, we will design an inexact semi-proximal augmented Lagrangian method (ALM) for the primal problem which attempts to exploit the block angular structure to solve the problem in parallel.
Our algorithm is motivated by the recent theoretical advances in
inexact semi-proximal  ALM  that is embedded in \cite{CST}.
In contrast to most existing augmented Lagrangian based decomposition algorithms
 where the solution for each subproblem
must be computed exactly or to very high accuracy, our algorithm
has the key advantage of allowing the subproblems to be solved approximately with progressive
accuracy. We will also elucidate the connection of our algorithm to the well-known
diagonal quadratic approximation (DQA)
algorithm of Mulvey and Ruszczy\'{n}ski \cite{MR92}.

In the pioneering work in \cite{KLM96}, an ADMM based framework was designed for the primal block angular problem
(P) wherein the variables are duplicated and auxiliary variables are
introduced to make the first ADMM subproblem in each iteration solvable in a distributed
fashion and that the succeeding second ADMM subproblem is a sufficiently simple quadratic program
which is assumed to be easy to solve.
However, the problem might still be difficult to solve if the scale of the original problem gets very large.
To overcome the potential computational inefficiency induced by the
extra variables and constraints, and also the relatively expensive step of having to solve a QP
subproblem in each iteration in the primal approach,
in this paper we will adopt the dual approach of solving (P). Specifically, we will design and implement a semi-proximal symmetric Gauss-Seidel based alternating direction method of multipliers (sGS-ADMM) to directly solve the dual problem, which will also solve the primal problem as a by-product.
The advantage of tackling the dual problem directly is that no extra variables are
introduced to decouple the constraints and no coupled QP subproblems are needed
to be solved in each iteration. We note that the sGS-ADMM
is an algorithm designed based on the recent advances in \cite{CST}; more details
will be presented later.

We consider the following primal block-angular optimization problem:
\begin{eqnarray*}
 \begin{array}{rl}
\mbox{(P}) \quad   \min & \sum_{i=0}^N f_i(x_i):= \theta_i(x_i) + \frac{1}{2}\inprod{x_i}{\cQ_i x_i} +
\inprod{c_i}{x_i} \\[8pt]
  \mbox{s.t.} &
 \underbrace{ \left[\begin{array}{ccccc}
 \cA_0 &   \cA_1 & \dots & \cA_N  \\[5pt]
  & \cD_1 &         &\vdots    \\[5pt]
       &    & \ddots &\vdots  \\[5pt]
       &    &       & \cD_N
  \end{array} \right]}_{\bB}
  \left[\begin{array}{c}
 x_0\\[5pt] x_1  \\[5pt] \vdots \\[5pt] x_N
  \end{array}\right]
  \;=\;
   \left[\begin{array}{c}
 b_0  \\[5pt] b_1 \\[5pt] \vdots \\[5pt] b_N
  \end{array}\right],
  \\[25pt]
   & x_i \in\cK_i, \; i=0,1,\ldots,N,
\end{array}
 \label{eq-1}
\end{eqnarray*}
where for each $i=0,1,\ldots,N$, $\theta_i : \cX_i \to (-\infty,\infty]$ is a proper closed convex function,
$\cQ_i : \cX_i \to \cX_i$ is a positive semidefinite linear operator, $\cA_i:\cX_i \to \cY_0$
and $\cD_i:\cX_i\to \cY_i$ are given linear maps, $c_i\in\cX_i$ and $b_i\in \cY_i$ are given data,
$\cK_i \subset \cX_i$ is a closed convex set that is typically a cone but not necessarily so,
and $\cX_i,\cY_i$ are real finite dimensional Euclidean spaces each equipped with an
inner product $\inprod{\cdot}{\cdot}$ and its induced norm $\norm{\cdot}$.
Note that the addition of the proper closed convex functions
in the objective gives us the flexibility to add nonsmooth terms such as $\ell_1$ regularization terms.
We should also mention that a constraint of the form
$b_i - \cD_i x_i \in  {\cal C}_i$, where ${\cal C}_i$ is a closed convex set can be put in
the form in (P) by introducing a slack variable $s_i$ so that
$[\cD_i, I] (x_i; s_i) = b_i$ and $(x_i;s_i)\in \cK_i\times {\cal C}_i$.

Without loss of generality, we assume that the constraint matrix $\bB$ in (P)
has full row-rank.
Let $n_i = {\rm dim}(\cX_i)$
and $m_i = {\rm dim}(\cY_i)$.
Observe that the problem (P) has $\sum_{i=0}^N m_i$ linear constraints and the dimension
of the decision variable is
$\sum_{i=0}^N n_i$. Thus even if $m_i$ and/or $n_i$ are moderate numbers, the
overall dimension of the problem can easily get very large when $N$ is large.

In the important special case of a block angular linear programming problem for which $\cQ_i=0$ and $\theta_i =0$ for all
$i=0,\ldots,N$,
the Dantzig-Wolfe decomposition method (which may be viewed as a
dual method based on the Lagrangian function
$\sum_{i=0}^N \inprod{c_i}{x_i}  - \inprod{u}{b_0- \sum_{i=0}^N \cA_i x_i}$) is a well known
classical approach for solving such a problem.
The Dantzig-Wolfe decomposition method has the attractive property that
in each iteration, $x_i$ can be computed individually from a smaller linear program (LP) for $i=1,\ldots,N$.
However, it is generally acknowledged that an augmented Lagrangian approach
has a number of important advantages over the usual Lagrangian  dual method.
For example, Ruszczy\'{n}ski stated in \cite{ruszcz95}  that the dual approach based on the ordinary Lagrangian can suffer from the nonuniqueness of the solutions of subproblems. In addition, solving the subproblem of the augmented Lagrangian approach would be  more
stable. In that paper, the well-known diagonal quadratic approximation (DQA) method is introduced. The DQA method is a very successful decomposition method and it has been a popular tool in stochastic programming. Thus it would be a worthwhile effort
to analyse it again to see whether further enhancements are possible.

To summarize, our first contribution  is in proposing several variants of augmented Lagrangian based algorithms for directly solving the primal form (P) of the convex composite quadratic conic programming problem with a primal block angular structure. We also show that they can be
considered as generalizations of the well-known DQA method. Our
 second contribution is in the  design and implementation of a specialized algorithm for solving the dual problem of (P). The algorithm is easy to implement and highly amenable to parallelization. Hence we expect it to be highly scalable for solving large scale problems with million of variables and constraints. Finally, we have proposed efficient implementations
 of the algorithms and conducted comprehensive numerical experiments to evaluate the
 performance of our algorithms against highly competitive state-of-the-art solvers
 in solving the problems (P) and (D).

This paper is organized as follows. We will derive the dual of the primal block angular problem (P)
 in section 2. In section 3, we will present our inexact semi-proximal augmented Lagrangian methods for the primal problem (P). In section 4, we will propose a semi-proximal symmetric Gauss-Seidel based ADMM for the dual problem of (P). For all algorithms we introduce, we conduct numerical experiments to evaluate their performance and the numerical results are reported in section 3.3 and  section 5. We conclude the paper in the final section.

\bigskip
\noindent {\bf Notation.}
\begin{itemize}
	\item We denote $[P;Q]$ or $(P;Q)$ as the matrix obtained by appending the matrix $Q$ to the last row of the matrix $P$, whereas we denote $[P,Q]$ or $(P,Q)$ as the matrix obtained by appending $Q$ to the last column of matrix $P$, assuming that they have the same number of columns or rows respectively. We also use the same notation  symbolically
	 for $P$ and $Q$ which are
	linear maps with compatible domains and co-domains.
		
	\item For any linear map  $\cT:\cX \to \cY$, we denote its adjoint as $\cT^*$. If $\cX = \cY$, and $\cT$ is self-adjoint and positive semidefinite, then for any $x \in \cX$ we have the notation $\|x\|_\cT:=\sqrt{\inprod{x}{\cT x}}$.
	
	\item Let $f:\cX \to (-\infty,+\infty]$ be an arbitrary closed proper convex function. We denote ${\rm dom} f$ as its effective domain and $\partial f$ as its subdifferential mapping.
The Fenchel conjugate function of $f$ is denoted as $f^*$.
	
	\item The Moreau-Yosida proximal mapping of $f$ is defined by ${\rm Prox}_f (y):=\arg\min_x \{f(x) + \frac{1}{2} \|x-y\|^2\}$.
\end{itemize}

%%%%%%%%%%%%%%%%%%%%%%%%%%%%%%%%%
\section{Derivation of the dual of (P)}

For notational convenience, we define
\begin{eqnarray}
\begin{array}{l}
 \bX = \cX_0\times \cX_1\times \cdots\times\cX_N, \quad
\bY = \cY_0\times \cY_1\times \cdots\times\cY_N,
\quad
\bK = \cK_0\times \cK_1\times \cdots\times \cK_N.
\end{array}
\label{eq-struct-1}
\end{eqnarray}
For each $x\in\bX$, $y\in\bY$, and $c\in\bX$, $b\in\bY$, we can express them as
\begin{eqnarray}
\begin{array}{l}
 x = (x_0; x_1; \cdots; x_N), \quad
y = (y_0; y_1; \cdots; y_N),
\\[5pt]
c = (c_0; c_1; \cdots; c_N), \quad
b = (b_0; b_1; \cdots; b_N).
\end{array}
\label{eq-struct-2}
\end{eqnarray}
We also define $\cA$, $\bQ$ and
$\mb{\theta}$ as follows:
\begin{eqnarray}
\begin{array}{c}
  \cA = [\cA_0,\cA_1,\ldots,\cA_N], \;
  \bQ(x) = \big(\cQ_0(x_0); \cQ_1(x_1); \ldots; \cQ_N(x_N)\big), \; \mb{\theta}(x) = \sum_{i=0}^N\theta_i(x_i).
\end{array}
\label{eq-struct-3}
\end{eqnarray}
Using the notation in \eqref{eq-struct-1}--\eqref{eq-struct-3}, we can write (P)
compactly in the form  of a general convex composite quadratic conic programming
problem:
\begin{eqnarray}
\begin{array}{ll}
 \min  \Big\{ \mb{\theta}(x) + \frac{1}{2}\inprod{x}{\bQ x} + \inprod{c}{x} \mid
\bB x-b =0, \;
x\in\bK \Big\}.
\end{array}
\label{eq-CCQP}
\end{eqnarray}
By introducing auxiliary variables $u,v \in \bX$, problem \eqref{eq-CCQP} can equivalently be written as the following model:
\begin{eqnarray}
\begin{array}{ll}
\min & \mb{\theta}(u) + \frac{1}{2}\inprod{x}{\bQ x} + \inprod{c}{x} + \delta_{\bK}(v)
\\[5pt]
\mbox{s.t.} & \bB x-b =0, \; u-x=0,\; v-x=0.
\end{array}
\label{eq-CCQP2}
\end{eqnarray}
To derive the dual of \eqref{eq-CCQP}, consider the following Lagrangian function
 for \eqref{eq-CCQP2}:
\begin{eqnarray*}
& & \cL (x,u,v; y,s,z) \\
&=& \mb{\theta}(u) + \frac{1}{2}\inprod{x}{\bQ x} + \inprod{c}{x} + \delta_{\bK}(v) - \inprod{y}{\bB x - b} - \inprod{s}{x - u} - \inprod{z}{x - v} \\
&=& \frac{1}{2}\inprod{x}{\bQ x} + \inprod{c - \bB^*y - s - z}{x} + \mb{\theta}(u) + \inprod{s}{u} + \delta_{\bK}(v) + \inprod{z}{v} + \inprod{y}{b},
\end{eqnarray*}
where $x, u, v, s, z \in \bX, y \in \bY$. Now
for a given subspace $\bW \subset \bX$ containing  ${\rm Range}(\bQ)$, the range space of $\bQ$, we have
\begin{eqnarray*}
\inf_{x} \cL (x,u,v; y,s,z)
&=& \inf_{x} \Big\{\frac{1}{2}\inprod{x}{\bQ x} + \inprod{c - \bB^*y - s - z}{x}\Big\} \\
&=&
\begin{cases}
	-\frac{1}{2}\inprod{w}{\bQ w}, \quad \text{if } c - \bB^*y - s - z
	= -\bQ w \text{ for some } w \in \bW, \\
	-\infty, \quad \text{otherwise}.
\end{cases}
\end{eqnarray*}
Also,
\begin{eqnarray*}
\inf_u \cL (x,u,v; y,s,z)
&=& \inf_u \big[ \mb{\theta}(u) + \inprod{s}{u}\big] = -\mb{\theta^*}(-s);\\
\inf_v \cL (x,u,v; y,s,z)
&=& \inf_v \big[\delta_{\bK}(v) + \inprod{z}{v}\big] = -\delta_{\bK}^*(-z).
\end{eqnarray*}
Hence the dual of \eqref{eq-CCQP2} is given by
\begin{eqnarray*}
& &\max_{y,s,z} \inf_{x,u,v}
\cL (x,u,v; y,s,z) \\
&=& \max_{y,s,z,w} \Big\{-\mb{\theta^*}(-s) -\frac{1}{2}\inprod{w}{\bQ w} + \inprod{y}{b} -\delta_{\bK}^*(-z) \mid \, -\bQ w + \bB^*y + s + z =c, \; w \in \bW
\Big\},
\end{eqnarray*}
or equivalently,
%By considering the Lagrangian function for the above problem, one can derive the following dual
%problem
%of \eqref{eq-CCQP}:
\begin{eqnarray}
 \begin{array}{rl}
-\min & \mb{\theta}^*(-s) + \frac{1}{2}\inprod{w}{\bQ w} - \inprod{b}{y} + \delta_{\bK}^*(-z)
\\[5pt]
\mbox{s.t.}& -\bQ w + \bB^*y + s+z = c,
\\[5pt]
& s\in \bX,\; y\in\bY,\; w \in\bW.
\end{array}
\label{eq-CCQP-dual}
\end{eqnarray}
It is not difficult to check that for all
$z =  (z_0;z_1;\ldots;z_N)$,
$s = (s_0;s_1;\ldots;s_N)\in \bX$, we have
\begin{eqnarray}
\begin{array}{l}
\delta_{\bK}^*(-z) =  \sum_{i=0}^N \delta_{\cK_i}^*(-z_i), \quad
 \mb{\theta}^*(-s) = \sum_{i=0}^N \theta_i^*(-s_i) .
\end{array}
\label{eq-struct-4}
\end{eqnarray}

Assume that both the primal and dual problems satisfy the (generalized) Slater's condition. Then the optimal solutions for both problems exist and they satisfy the following Karush-Kuhn-Tucker (KKT) optimality conditions:
\begin{eqnarray}
\begin{cases}
\begin{array}{ll}
\bB x-b =0, &  \\
-\bQ w + \bB^*y + s+z - c= 0, \quad \bQ w - \bQ x = 0, \quad w \in \bW,   \\
-s \in \partial  \mb{\theta}(x) \;\Leftrightarrow\; x - {\rm Prox}_{\mb{\theta}} (x-s) = 0, &\\
x-\Pi_{\bK}(x-z)=0.  &
\end{array}
\end{cases}
\label{KKT}
\end{eqnarray}
By applying the structures in \eqref{eq-struct-1}--\eqref{eq-struct-3} and
\eqref{eq-struct-4} to
\eqref{eq-CCQP-dual}, we get explicitly the dual of (P):
\begin{eqnarray}
\begin{array}{ll}
\mbox{(D)} \quad -\min &  \sum_{i=0}^N \theta_i^*(-s_i) + \delta_{\cK_i}^*(-z_i) + \frac{1}{2}\inprod{w_i}{\cQ_i (w_i)} -\inprod{b_i}{y_i}
\\[10pt]
\mbox{s.t.} &
\left[\begin{array}{c} \cA_0^* \\[5pt] \cA_1^* \\[5pt] \vdots \\[5pt] \cA_N^* \end{array}\right] y_0 +
\left[ \begin{array}{l}
\qquad\;\;  -\cQ_0 w_0 + s_0 + z_0 \\[5pt]
\cD_1^* y_1 -\cQ_1 w_1 + s_1+ z_1\\[5pt]
\qquad\qquad \vdots\qquad\qquad \\[5pt] \cD_N^* y_N -\cQ_N w_N + s_N + z_N
\end{array}\right] = c,
\\[40pt]
& w_i \in \cW_i,\; i=0,1,\ldots,N,
\end{array}
\label{eq-struct-CCQP-dual}
\end{eqnarray}
where for each $i=0,1,\ldots,N$, $\cW_i\subset \cX_i$ is a given subspace containing ${\rm Range}(\cQ_i)$.

%%%%%%%%%%%%%%%%%%%%%%%%%%%%%%%%%%%
\section{Inexact semi-proximal augmented Lagrangian methods for the primal problem
(P)}

First we rewrite (P) in the following form:
\begin{eqnarray}
 \min\Big\{\sum_{i=0}^N f_i(x_i) + \delta_{F_i}(x_i)  \;\mid\; \cA x = b_0,\; x= (x_0;x_1;\ldots;x_N) \in \bX
 \Big\},
 \label{eq-P2}
\end{eqnarray}
where
$F_0 = \cK_0$, and
$F_i = \{ x_i\in\cX_i\mid \cD_i x_i  = b_i, x_i\in \cK_i\}$, $i=1,\ldots,N$.
For a given parameter $\sig > 0$, we consider
the following augmented Lagrangian function associated with (\ref{eq-P2}):
\begin{eqnarray}
\begin{array}{lll}
  L_\sig(x;y_0) = \sum_{i=0}^N f_i (x_i) + \delta_{F_i}(x_i) +
  \frac{\sig}{2} \norm{\cA x - b_0 -\sig^{-1}y_0 }^2 -\frac{1}{2\sig}\norm{y_0}^2.
\end{array}
  \label{eq-1-L}
\end{eqnarray}
The augmented Lagrangian method for solving (\ref{eq-P2}) has the following template.
\\[5pt]
\centerline{\fbox{\parbox{\textwidth}{
{\bf ALM}.
Given $\sig > 0$ and $y_0^0$
$\in \mathcal{Y}_0$.
Perform the following steps in each iteration.
\begin{description}
\item[Step 1.]
$
 x^{k+1} \approx \mbox{argmin}_x \; L_\sig(x;y_0^k) .
$
\item[Step 2.] $y_0^{k+1} = y_0^k + \tau\sig(b_0- \cA x^{k+1})$,
where $\tau\in (0,2)$ is the step-length.
\end{description}
}}}

\bigskip
As one may observe from Step 1 of the ALM,
an undesirable feature  in the method is that it
destroys the separable structure in the Dantzig-Wolfe decomposition method.
Although the feasible sets for the $x_i$'s are separable, the objective function
has a quadratic term which couples all the $x_i$'s.

Here we propose to add a semi-proximal term to the augmented
Lagrangian function to overcome the difficulty of non-separability. In this case, the function $L_\sig(x; y_0^k)$ in Step 1 of the ALM
is majorized by an additional semi-proximal term at the point $x^k$, i.e.,
\begin{eqnarray*}
  L_\sig(x; y_0^k)  + \frac{\sig}{2}\norm{x-x^k}_{\cT}^2,
\end{eqnarray*}
where $\cT$ is a given positive semidefinite self-adjoint linear operator which should be
chosen appropriately to decompose the computation of the $x_i$'s in Step 1 of  the ALM
while at the same time the added proximal term should be as small as possible.
In this paper, we choose $\cT$ to be the following positive semidefinite linear operator:
\begin{eqnarray}
  \cT = \diag{\cJ_0, \ldots, \cJ_N}
   - \cA^*\cA,
\label{eq-cT}
\end{eqnarray}
where $\cJ_i \succeq \beta_i I + \cA_i^*\cA_i$, with $\beta_i = \sum_{j=0, j\not = i}^N \norm{\cA_i^* \cA_j}_2$
for each $i=0,1,\ldots, N.$
Such a choice is generally less conservative than the usual choice
of $\widehat{\cT}$ given in \eqref{eq-cT-2}. It is especially
a good choice  when $\cA_i$ and $\cA_j$ are nearly
orthogonal to one another for most of the index pairs $(i,j)$.

With the choice in \eqref{eq-cT}, we get
\begin{eqnarray*}
 && L_\sig(x; y_0^k)  + \frac{\sig}{2}\norm{x-x^k}_{\cT}^2 \\
 &=& \sum_{i=0}^N \big(f_i (x_i) + \delta_{F_i}(x_i) \big) +
 \frac{\sig}{2} \norm{\cA x - b_0 -\sig^{-1}y_0 }^2 -\frac{1}{2\sig}\norm{y_0}^2  + \frac{\sig}{2}\norm{x-x^k}_{\cT}^2 \\
% &=& \sum_{i=0}^N \big(f_i (x_i) + \delta_{F_i}(x_i)\big) -\frac{1}{2\sig}\norm{y_0}^2 \\
% && + \frac{\sig}{2}\Big[ \inprod{x}{(\cA^* \cA + \cT)x} - 2\inprod{x}{\cA^* (b_0 + \sig^{-1}y_0) + \cT x^k} + \norm{b_0 +\sig^{-1}y_0}^2 + \norm{x^k}_\cT \Big]  \\
 &=& \sum_{i=0}^N \Big(f_i (x_i) + \delta_{F_i}(x_i) + \frac{\sig}{2}\Big[\inprod{x_i}{\cJ_i x_i} - 2\inprod{x_i}{\cA_i^*(b_0 + \sig^{-1}y_0 - \cA x^k) +
 	\cJ_i x_i^k}\Big]\Big)
 	\\
&&+ \frac{\sig}{2}\Big[\norm{b_0 +\sig^{-1}y_0}^2 + \norm{x^k}_\cT\Big]
 -\frac{1}{2\sig}\norm{y_0}^2.
\end{eqnarray*}

The inexact semi-proximal ALM (sPALM) we consider for solving the
primal block angular problem (P)  through \eqref{eq-P2} is given as follows.
\\[5pt]
\centerline{\fbox{\parbox{\textwidth}{
{\bf sPALM.} Given $\sig >0$ and $y_0^0 \in \mathcal{Y}_0$.
Let $\{\varepsilon_k\}$ be a given
summable sequence of nonnegative numbers.
Perform the following steps in each iteration.
\begin{description}
\item[Step 1.] Compute
\begin{eqnarray}
x^{k+1}\; \approx\; \widehat{x}^{k+1}  := \mbox{argmin}\{ L_\sig(x;y_0^k) + \frac{\sig}{2}\norm{x-x^k}^2_\cT \},
\label{eq-sPALM-1}
\end{eqnarray} with residual
\begin{eqnarray}
 d^{k+1} \in \partial_x L_\sig(x^{k+1}; y_0^k) + \sig \cT(x^{k+1}-x^k), \quad
 \label{eq-stop}
\end{eqnarray}
satisfying $\norm{d^{k+1}} \leq \varepsilon_k$.
Let
 $\cG_i = \cQ_i+\sig\cJ_i$,
$g^k_i = \cQ_i x_i^k + c_i +\sig  \cA_i^*(\cA x^k-b_0 -\sig^{-1}y_0^k) - \cG_i x_i^k$.
Due to the separability of the variables in \eqref{eq-sPALM-1} because of the specially chosen $\cT$, one can compute {\bf in parallel} for  $i=0,1,\ldots,N$,
\begin{eqnarray}
 x^{k+1}_i \;\approx\; \widehat{x}^{k+1}_i  :=  \mbox{argmin} \Big\{  \theta_i(x_i)+
\frac{1}{2}\inprod{x_i}{\cG_i x_i} + \inprod{g_i^k}{x_i}
\,\mid\, x_i\in F_i\Big\},
\label{eq-sPALM-sub}
\end{eqnarray}
with the residual $d_i^{k+1} := v_i^{k+1}+ \cG_i x^{k+1}_i + g_i^k$ for
some $v^{k+1}_i \in \partial (\theta_i+\delta_{F_i})(x_i^{k+1}) $ and satisfying
\begin{eqnarray}
\norm{d_i^{k+1}} \;\leq\;  \frac{1}{\sqrt{N+1}} \varepsilon_k.
\end{eqnarray}
\item[Step 2.] $y_0^{k+1} = y_0^k + \tau\sig(b_0-\cA x^{k+1})$, where
$\tau\in (0,2)$ is the steplength.
\end{description}
}}}

\bigskip
Observe that  with the introduction of the
semi-proximal term $\frac{\sig}{2}\norm{x-x^k}_\cT^2$ to the augmented Lagrangian
function
in Step 1 of the sPALM,
we have decomposed the large coupled problem involving
$x$ in ALM into $N+1$ smaller independent problems that can be
solved in parallel. For the case of a quadratic or linear program,
we can employ a  powerful solver such as   Gurobi or Mosek to
efficiently solve these smaller problems.

In order to judge how accurately the decomposed subproblems
in Step 1 must be solved, we need to analyse the stopping condition for \eqref{eq-sPALM-sub} in
detail.
In particular, we need to
find $v_i^{k+1} \in \partial (\theta_i+\delta_{F_i})(x_i^{k+1})$ for $i=0,1,\ldots,N$.
This can be done by considering
 the dual of the subproblem \eqref{eq-sPALM-sub}, which could be written as:
\begin{eqnarray}
\begin{array}{rl}
-\min & \theta_i^*(-s_i) + \frac{1}{2}\inprod{w_i}{\cG_i w_i} - \inprod{b_i}{y_i} + \delta_{\cK_i}^*(-z_i)
\\[5pt]
\mbox{s.t.}& -\cG_i w_i + \cD_i^* y_i + s_i + z_i = g_i^k,
\\[5pt]
& s_i \in \cX_i,\;\; y_i\in\cY_i,\; w_i \in\cW_i, \quad i=1,\ldots,N.
\end{array}
\label{eq-sPALM-sub-dual}
\end{eqnarray}
Note that for $i=0$, we have a similar problem as the above but the terms involving $y_i$ are absent. For the discussion below, we will just focus on the case where $i=1,\ldots,N$, the case
for $i=0$ can be derived similarly.
One can estimate $v_i^{k+1}$ to be $-\cD_i^* y_i^{k+1} - s_i^{k+1}
- z_i^{k+1}$ for a computed dual solution $(y_i^{k+1},s_i^{k+1},z_i^{k+1})$
and the residual $d_i^{k+1}$ is simply the
residual in the dual feasibility constraint in the above problem.

\begin{remark}
In the sPALM, some of the dual variables for (D) are not explicitly constructed. Here we describe how
they can be estimated. Recall that for (D), we want to get
\[
-\cQ_i x_i + \cA_i^*y_0 + \cD_i^* y_i + s_i + z_i - c_i =0\quad \forall\; i=0,1,\ldots,N.
\]
Note that for convenience, we introduced $\cD_0^* =0$.
From the KKT conditions for \eqref{eq-sPALM-sub} and \eqref{eq-sPALM-sub-dual}, we have that
\begin{eqnarray*}
&-\cG_i w_i^{k+1} + \cD^*_i y_i^{k+1} + s_i^{k+1} + z_i^{k+1}  - g_i^k
 =: R^d_i \;\approx\; 0,
 &
 \\[5pt]
 &\cG_i w_i^{k+1} - \cG_i x_i^{k+1} \; \approx\;  0.&
\end{eqnarray*}
By using the expression for $\cG_i$, $g_i^k$ and $y_0^{k+1}$, we get
\begin{eqnarray*}
&& \hspace{-0.7cm}
-\cQ_i x_i^{k+1} + \cA_i^* y_0^{k+1}
 + \cD_i^* y_i^{k+1} + s^{k+1}_i + z^{k+1}_i -c_i
 \\[5pt]
&=& R^d_i + (\cG_i w_i^{k+1} - \cG_ix_i^{k+1}) + \sig \cJ_i (x_i^{k+1}-x_i^{k})
+ \sig\cA_i^* \cA(x^k - x^{k+1})
+  (\tau-1)\sig \cA_i^*( b_0-\cA x^{k+1}).
\end{eqnarray*}
Note that the right-hand-side quantity in the above equation will converge to $0$ based
on the convergence of sPALM and the KKT conditions for \eqref{eq-sPALM-sub} and
\eqref{eq-sPALM-sub-dual}.
Thus by using the dual variables computed from solving \eqref{eq-sPALM-sub-dual}, we can
generate the dual variables for (D).

\end{remark}

%%%%%%%%%%%%%%%%%%%%%%%%%%%%%%%%
\subsection{Convergence of the inexact sPALM}

The convergence of the inexact sPALM for solving \eqref{eq-P2}
can be established readily by using known
results in \cite{CLST}. To do that, we need to first reformulate \eqref{eq-P2}
into the form required in \cite{CLST}
as follows:
\begin{eqnarray}
 \min\Big\{h(x) + \psi(x)   \;\mid\; \cA x = b_0,\; x= (x_0;x_1;\ldots;x_N) \in \bX
 \Big\},
 \label{eq-P22}
\end{eqnarray}
where $h(x) = \sum_{i=0}^N \frac{1}{2}\inprod{x_i}{Q_ix_i} + \inprod{c_i}{x_i}$
and $\psi(x) = \sum_{i=0}^N \theta_i(x_i) + \delta_{F_i}(x_i)$.
Its corresponding KKT residual mapping is given by
\begin{eqnarray}
  {\cal R} (x,y_0) = \left( \begin{array}{c}
   b_0 - \cA x \\[3pt]
   x - {\rm Prox}_\psi (x - \bQ x - c - \cA^* y_0)
  \end{array}\right)
  \quad \forall\; x\in \bX, \; y_0 \in \cY_0.
\end{eqnarray}
Note that $(x,y_0)$ is a solution of the KKT system of \eqref{eq-P22} if and only
${\cal R}(x,y_0) = 0$.

Now we state the global convergence theorem here
for the convenience of the readers.
Define the self-adjoint positive definite linear operator $\cV:\bX\to\bX$ by
$$
\cV:=\tau\sigma\Big(\bQ+ \sig\cT+\frac{2-\tau}{6}\sigma\cA\cA^*\Big).
$$
We have the following convergence result for the inexact sPALM.
\begin{theorem}
\label{thm:main}
Assume that the solution set to the KKT system of \eqref{eq-P2} is nonempty
and $(\overline x, \overline{y}_0)$ is  a solution. Then, the sequence $\{(x^k,y_0^k)\}$ generated by sPALM is well-defined such that for any $k\ge 1$,
\[\label{estd}
\|x^{k+1}-\widehat{x}^{k+1}\|^2_{\bQ+\sig\cT+\sigma\cA\cA^*}
\le
\langle d^{k+1}, x^{k+1}-\widehat{x}^{k+1}\rangle,
\]
and for all $k=0,1,\ldots$,
 \begin{eqnarray*}
\label{ineq:keyinequality}
\begin{array}{l}
\left(\| x^{k+1} - \overline{x}\|_{\widehat{\cV}}^2
+\norm{y_0^{k+1} -\overline{y}_0}^2
\right)
-\left( \| x^{k} - \overline{x}\|_{\widehat{\cV}}^2
+ \norm{y_0^{k} -\overline{y}_0}^{2}
\right)
\\[8pt]
\le
-\left( \frac{2-\tau}{3\tau}\|y_0^{k}-y_0^{k+1}\|^{2}
+\|x^{k+1}-x^{k}\|_{\cV}^2
-2\tau\sigma\langle d^k,x^{k+1}-\overline{x} \rangle\right),
\end{array}
\end{eqnarray*}
where $\widehat{\cal V} = {\cal V} + \frac{2-\tau}{6}\sig\cA\cA^*$.
Moreover,  the sequence $\{(x^k,y_0^k)\}$
converges to a solution to the KKT system of \eqref{eq-P2}.
\end{theorem}
\begin{Proof}
	The result can be proved directly from the convergence result in \cite[Theorem 1]{CLST}.
\end{Proof}

\medskip
The local linear convergence of sPALM can also be established if the KKT
residual mapping ${\cal R}$ satisfies the following error bound condition:
there exist positive constants $\kappa$ and $r$  such that
${\rm dist}((x,y_0),\Omega) \leq \kappa \norm{{\cal R}(x,y_0)}$
for all $(x,y_0)$ satisfying $\norm{(x,y_0)-(x^*,y_0^*)} \leq r$,
where $\Omega$ is the solution set of \eqref{eq-P22} and $(x^*,y_0^*)$ is
a particular solution of \eqref{eq-P22}.
In order to save some space, we will not state the theorem here but
refer the reader to \cite[Theorem 2]{CLST}.

%%%%%%%%%%%%%%%%%%%%%
\subsection{Comparison of sPALM with the diagonal quadratic approximation method
and its recent variants}

Let  $\rho := (N+1)^{-1}$.
Consider the following linear operator
\begin{eqnarray}
 \widehat{\cT} = \diag{\cE_0,\ldots,\cE_N} - \cA^*\cA,
\label{eq-cT-2}
\end{eqnarray}
where $\cE_i \succeq \rho^{-1}\cA_i^*\cA_i$ for all $i=0,1,\ldots,N$.
It is not difficult to show
that $\widehat{\cT} \succeq 0$.
%Consider the convex function $\phi(z) = \norm{z}^2$ where $z=(z_0,z_1,\ldots,z_N)\in \bY$. It is clear from the convexity of $\phi$ that
%\begin{eqnarray*}
%\Big\|\mbox{$\sum_{i=0}^N$} \rho z_{i}\Big\|^2 =  \phi\Big(\mbox{$\sum_{i=0}^N$} \rho z_{i}\Big)
%\leq \mbox{$\sum_{i=0}^N$} \rho\phi(z_i) = \mbox{$\sum_{i=0}^N$} \rho \norm{z_i}^2.
%\end{eqnarray*}
%Hence  $\norm{\sum_{i=0}^N z_i}^2 \leq \rho^{-1} \sum_{i=0}^N \norm{z_i}^2$.
%By taking $z_i = \cA_i x_i$, we get
%\begin{eqnarray*}
%  \inprod{x} {\cA^*\cA x}= \norm{\mbox{$\sum_{i=0}^N$} \cA_i x_i}^2 \leq \rho^{-1} \mbox{$\sum_{i=0}^N$} \norm{\cA_ix_i}^2 = \rho^{-1}\mbox{$\sum_{i=0}^N$} \inprod{x_i}{\cA_i^*\cA_i x_i},
%\end{eqnarray*}
%which implies that $\cA^*\cA \preceq  \mbox{diag}(\rho^{-1}\cA_0^* \cA_0, \cdots, \rho^{-1}\cA_N^* \cA_N) \preceq \diag{\cE_0,\ldots,\cE_N}$.

If instead of \eqref{eq-cT}, we choose $\cT$ to be the linear operator  given
in \eqref{eq-cT-2}, then instead of {\bf sPALM},
we get the following variant of the inexact sPALM.
\\[5pt]
\centerline{\fbox{\parbox{\textwidth}{
{\bf sPALM-b.}
 Given $\sig >0$ and $y_0^0$
 $\in \mathcal{Y}_0$.
Let $\{\varepsilon_k\}$ be a given
summable sequence of nonnegative numbers.
Perform the following steps in each iteration.
\begin{description}
\item[Step 1.]  Let
$g^k_i = \cQ_i x_i^k+c_i +\sig  \cA_i^*(\cA x^k-b_0 -\sig^{-1}y_0^k)
- (\cQ_i+\sig\cE_i) x_i^k$. Compute (in parallel) for  $i=0,1,\ldots,N$,
\begin{eqnarray}
 x^{k+1}_i \;\approx\; \mbox{argmin} \Big\{  \theta_i(x_i)+
\frac{1}{2}\inprod{x_i}{(\cQ_i +\sig\cE_i) x_i} + \inprod{g_i^k}{x_i}
\,\mid\, x_i\in F_i\Big\},
\label{eq-sPALM-b}
\end{eqnarray}
with the residual $d_i^{k+1} := v_i^{k+1} + (\cQ_i + \sig\cE_i) x^{k+1}_i + g_i^k$ for
some $v_i^{k+1} \in \partial (\theta_i + \delta_{F_i})(x^{k+1}_i)$ and satisfying
$\norm{d_i^{k+1}} \leq \frac{1}{\sqrt{N+1}}\varepsilon_k.$
\item[Step 2.] $y_0^{k+1} = y_0^k + \tau\sig(b_0-\cA x^{k+1})$, where
$\tau\in (0,2)$ is the steplength.
\end{description}
}}}

\bigskip
In \cite{Rusz89}, Ruszczy\'{n}ski
proposed the diagonal quadratic approximation (DQA) augmented Lagrangian
method
that aims to
solve a problem of the form (P).
As already mentioned, the DQA method is a very successful decomposition method
that is frequently used in stochastic programming.
%To describe the DQA method, it is
%more convenient to write (P) in the following form:
%$$
%\begin{array}{l}
% \min \Big\{ \sum_{i=0}^N f_i(x_i) \mid  \cA x-b_0=0, x_i\in F_i, \; i=0,1,\ldots,N\Big\}.
%\end{array}
%$$
Although it was not derived in our way in \cite{Rusz89},
we shall see later that the DQA method  can roughly be derived as the augmented Lagrangian method described in {\bf ALM}
where the minimization problem in Step 1
is solved approximately by a proximal gradient method, with the proximal term
chosen specially using the linear operator $\widehat{\cT}$ in
\eqref{eq-cT-2}  to make the resulting subproblem separable.
\\[5pt]
\centerline{\fbox{\parbox{\textwidth}{
{\bf ALM-DQA-mod.}
Given $\sig >0$, $y_0^0\in \mathcal{Y}_0$
and $x^0\in \bX$.
Let $\{\varepsilon_k\}$ be a given
summable sequence of nonnegative numbers.
Perform the following steps in each iteration.
\begin{description}
\item[Step 1.]
Starting with $\hat{x}^0 = x^k$, iterate the following step for $s=0,1,\ldots$ until convergence:
\begin{itemize}
\item Compute $
    \hat{x}^{s+1} \approx \mbox{argmin} \big\{  L_\sig(x; y_0^k)
 + \frac{\sig}{2} \norm{x-\hat{x}^s}^2_{\widehat{\cT}} \mid x\in \bX\big\}.
$
As the problem is separable, one can compute in parallel for $i=0,1,\ldots,N$,
\begin{eqnarray}
    \hat{x}_i^{s+1}
&\approx & \mbox{argmin} \big\{ f_i(x_i) +
 \frac{\sig}{2}\inprod{x_i-\hat{x}_i^s}{\cE_i(x_i-\hat{x}_i^s)}
+ \inprod{x_i-\hat{x}_i^s}{\hat{g}_i^s} \mid x_i\in F_i
  \big\}
\nn \\[5pt]
&=&\mbox{argmin} \big\{ \theta_i(x_i) +
 \frac{1}{2}\inprod{x_i}{(\cQ_i+\sig\cE_i) x_i}
+ \inprod{x_i}{\bar{g}_i^s} \mid x_i\in F_i
  \big\},
\label{eq-DQA-2}
\end{eqnarray}
where $\hat{g}^s_i =\sig \cA^*_i(\cA\hat{x}^s-b_0-\sig^{-1}y_0^k)$,
$\bar{g}^s_i =  \cQ_i \hat{x}^s_i + c_i+\sig \cA_i^*(\cA\hat{x}^s-b_0-\sig^{-1}y_0^k) - (\cQ_i+\sig\cE_i)\hat{x}_i^s$.
\end{itemize}
At termination, set $x^{k+1} = \hat{x}^{s+1}$.
\item[Step 2.] $y_0^{k+1} = y_0^k + \tau\sig(b_0-\cA x^{k+1})$, where
$\tau\in (0,2)$ is the steplength.
\end{description}
}}}

\bigskip
Observe that the subproblem \eqref{eq-sPALM-b} in Step 1 of sPALM-b is
exactly one step of the proximal gradient method \eqref{eq-DQA-2}
in Step 1 of the ALM-DQA-mod. As solving the problem of the form in \eqref{eq-DQA-2}
multiple times for each iteration of the
ALM-DQA-mod may be expensive, it is highly conceivable that the overall efficiency of
sPALM-b could be better than that of the ALM-DQA-mod.

Next, we elucidate the connection between {\bf ALM-DQA-mod} and the DQA method
described in \cite{Rusz89}.
Given $\hat{x}_i^s\in F_i$, we can parameterize a given $x_i$  as
$$
 x_i = \hat{x}_i^s + \rho d_i = (1-\rho) \hat{x}_i^s + \rho (\hat{x}_i^s + d_i), \;
 i=0,1,\ldots,N,
$$
with $\rho=(N+1)^{-1}\in (0,1]$.
Then by convexity, $f_i(x_i)  \leq (1-\rho)f _i(\hat{x}_i^s) + \rho f_i(\hat{x}_i^s + d_i)$.
Also,  if $\hat{x}_i^s + d_i \in F_i$, then $x_i \in F_i$ since $\hat{x}_i^s\in F_i$.
From here, we have that for all $x\in F_0 \times F_1\times \cdots \times F_N$,
\begin{eqnarray}
&& \hspace{-0.7cm} L_\sig(x; y_0^k)
+ \frac{\sig}{2} \norm{x-\hat{x}^s}^2_{\widehat{\cT}} + \frac{1}{2\sig}\norm{y_0^k}^2
\nn \\[5pt]
&\leq &(1-\rho) \sum_{i=0}^N f_i(\hat{x}_i^s)
+ \rho \sum_{i=0}^N f_i(\hat{x}_i^s + d_i) + \frac{\sig}{2}\norm{\cA (\hat{x}^s + \rho d) -b_0-\sig^{-1} y_0^k}^2  + \frac{\sig\rho^2}{2} \norm{d}^2_{\widehat{\cT}}
\nn \\[5pt]
&=&  \sum_{i=0}^N\rho f_i(\hat{x}_i^s + d_i) + \rho\inprod{d_i}{\hat{g}_i^s}
+ \frac{\sig\rho^2}{2}\inprod{d_i}{\cE_i d_i}
  +(1-\rho) \sum_{i=0}^N f_i(\hat{x}_i^s)
+\frac{\sig}{2}\norm{\cA\hat{x}^s - b_0 - \sig^{-1} y_0^k}^2.\qquad
\label{eq-L-major}
\end{eqnarray}
Hence instead of
\eqref{eq-DQA-2}, we may consider to minimize the majorization of
$L_\sig(x;y^k_0) + \frac{\sig}{2}\norm{x-\hat{x}^s}_{\widehat{\cT}}^2$ in \eqref{eq-L-major}, and
compute for  $ i=0,1,\ldots,N$,
\begin{eqnarray}
    d_i^{s+1}
=  \mbox{argmin} \, \rho \big\{  f_i(\hat{x}_i^s  + d_i) +
 \frac{{\sig \rho}}{2}\inprod{d_i}{\rho\cE_i d_i}
+ \inprod{d_i}{\hat{g}_i^s} \mid \hat{x}^s_i + d_i \in F_i, \; d_i\in\cX_i
  \big\}.\quad
\label{eq-DQA-d}
\end{eqnarray}
We get the DQA method
of \cite{Rusz89}
if we
 take $\cE_i =\rho^{-1}\cA_i^*\cA_i$,
 compute $d^{s+1}$ exactly  in the above subproblem \eqref{eq-DQA-d}, and set
$$
  \hat{x}^{s+1}_ i = \hat{x}^s_i + \rho d_i^{s+1}, \quad i = 0,1,...,N,
$$
instead of the solution in \eqref{eq-DQA-2}.
Thus we may view the DQA method as an augmented Lagrangian method
for which the subproblem in Step 1 is solved by a majorized proximal gradient method
with the proximal term chosen to be $\frac{\sig}{2}\norm{x-\hat{x}^s}^2_{\widehat{\cT}}$
in each step.

\begin{remark} When the $\cA_i$'s are matrices, the majorization $\cA^*\cA \preceq \diag{\cE_0,\ldots,\cE_N}$ can be improved as follows, as has been done in \cite{CDZ}.
Let
$$
 I_j = \big\{ i\in\{0,1,\ldots,N\} \mid e_j^T \cA_i \not = 0 \big\}, \quad
 \chi := \max\{ |I_j|\mid j=1,\ldots,m\} \;\leq\; N+1 = \rho^{-1}.
$$
Then
\begin{eqnarray*}
\begin{array}{l}
 \norm{\cA x}^2 = \norm{\sum_{i=0}^N \cA_i x_i}^2
= \sum_{j=1}^m
 |\sum_{i=0}^N e_j^T\cA_i x_i|^2
= \sum_{j=1}^m
\abs{\sum_{i\in I_j} e_j^T\cA_i x_i}^2
\\[8pt]
\;\leq\; \sum_{j=1}^m \Big( |I_j|\sum_{i\in I_j} \abs{e_j^T\cA_i x_i}^2 \Big)
\;\leq\;
\chi \sum_{j=1}^m \sum_{i\in I_j} \abs{e_j^T\cA_i x_i}^2
\\[8pt]
= \chi \sum_{j=1}^m \sum_{i=0}^N \abs{e_j^T\cA_i x_i}^2
= \chi  \sum_{i=0}^N \norm{\cA_i x_i}^2.
\end{array}
\end{eqnarray*}
That is, $\cA^*\cA \preceq\mbox{diag}( \chi \cA_0^* \cA_0, \cdots,  \chi \cA_N^* \cA_N)$. Such an improvement has been considered in
\cite{CDZ}. It is straightforward to incorporate the improvement into ALM-DQA-mod
by simply replacing
$\cE_i = \rho^{-1}\cA_i^*\cA_i$ in \eqref{eq-cT-2} by $\chi \cA_i^*\cA_i$ for each $i=0,1,\ldots,N.$
\end{remark}

With the derivation of the {\bf ALM-DQA-mod} as an augmented Lagrangian method with its
subproblems solved by a specially chosen proximal gradient method, we can
leverage on this viewpoint to
design an accelerated variant of this method. Specifically, we can improve the
efficiency in solving the subproblems by using an inexact accelerated
proximal gradient (iAPG) method, and we will also use a proximal term
based on the linear operator \eqref{eq-cT}, which
is typically less conservative than the term $\frac{\sig}{2}\norm{x-x^k}_{\widehat{\cT}}$
used in the DQA method.
\\[5pt]
\centerline{\fbox{\parbox{\textwidth}{
{\bf ALM-iAPG.}
Given $\sig >0$, $y_0^0\in \mathcal{Y}_0$ and $x^0\in \bX$.
Let $\{\varepsilon_k\}$ be a given
summable sequence of nonnegative numbers.
Perform the following steps in each iteration.
\begin{description}
\item[Step 1.]
Starting with $\hat{x}^0 = \bar{x}^0= x^k$, $t_0=1$, iterate the following step for $s=0,1,\ldots$ until convergence:
\begin{itemize}
\item Compute
  $ \hat{x}^{s+1} \approx \mbox{argmin} \big\{  L_\sig(x; y_0^k)
 + \frac{\sig}{2} \norm{x-\bar{x}^s}^2_{\cT} \mid x_i\in F_i, i=0,1,\ldots,N\big\}.$
As the problem is separable, one can compute in parallel for $i=0,1,\ldots,N$,
\begin{eqnarray}
\begin{array}{lll}
    \hat{x}_i^{s+1}
&\approx \mbox{argmin} \Big\{ \theta_i(x_i) +
 \frac{1}{2}\inprod{x_i}{\cG_i x_i}
+ \inprod{x_i}{\bar{g}_i^s} \mid x_i\in F_i
  \Big\},
\end{array}
\label{eq-DQA-3}
\end{eqnarray}
where $\cG_i = \cQ_i+\sig\cJ_i$,
$\bar{g}^s_i = \cQ_i \bar{x}^s_i + c_i +\sig  \cA_i^*(\cA \bar{x}_i^s-b_0 -\sig^{-1}y_0^k) - \cG_i \bar{x}_i^s$.
\item Compute $t_{s+1} = (1+\sqrt{1+4 t_s^2})/2$, $\beta_{s+1} = (t_s-1)/t_{s+1}$.
\item Compute $\bar{x}^{s+1} = (1+\beta_{s+1})\hat{x}^{s+1} - \beta_{s+1} \hat{x}^s$.
\end{itemize}
At termination, set $x^{k+1} = \hat{x}^{s+1}$.
\item[Step 2.] $y_0^{k+1} = y_0^k + \tau\sig(b_0-\cA x^{k+1})$, where
$\tau\in (0,2)$ is the steplength.
\end{description}
}}}

%%%%%%%%%%%%%%%
\subsection{Numerical performance of sPALM and ALM-DQA-mod}

In this subsection, we compare the performance of the sPALM and ALM-DQA-mod algorithms for solving several linear and quadratic test instances. The detailed description of the datasets
is given in Section \ref{sec-5}.
We also report the number of constraints and variables of the instances in the table. For all the instances, we have $m_1=m_2=...=m_N$ and $n_1=n_2=...=n_N$. Hence we denote them as $m_i$ and $n_i$ respectively.

Table \ref{tablePrimal} compares the performance of the two solvers sPALM and ALM-DQA-mod for the primal problem (P) through \eqref{eq-P2} against that of the solver sGS-ADMM for the dual problem \eqref{eq-struct-CCQP-dual}.
 The details of the dual approach will be presented in the next section. Here, we could observe that sPALM and ALM-DQA-mod always require much
longer runtime to achieve the same accuracy level in the
relative KKT residual when compare to sGS-ADMM, although the former algorithms generally take a smaller number of outer iterations. In addition, the ALM-DQA-mod algorithm is slightly slower than sPALM on the whole though the difference is not
too significant.
%On the other hand, ALM-DQA-mod is 2-3 times faster than sPALM for the last two larger instances.
Note that our preliminary implementation of the algorithms is in {\sc Matlab} which does not have a good
support for parallel computing. In a full scale implementation, one may try to implement these algorithms on an appropriate parallel computing platform with a good parallelization support. Nevertheless, the
inferior performance of the two primal approaches has motivated us to instead
consider the dual approach of
designing an efficient algorithm for the dual problem \eqref{eq-struct-CCQP-dual}.

\setlength\tabcolsep{2.6pt}
\begin{center}
	\scriptsize
	\begin{longtable}{| l | c | c | c || c | c || c | c || c | c|}	
						\caption{\footnotesize
			Comparison of computational results between sGS-ADMM and two variants of ALM for primal block angular problem. All the run result are obtained using \textbf{single thread}. Here, ``Iter'' is the number of outer iterations performed, and ``Time'' is the total runtime in seconds.
			}
		\label{tablePrimal}

		\\
		\hline
		\multicolumn{4}{|c||}{} & \multicolumn{2}{c||}{sGS-ADMM} &\multicolumn{2}{c||}{sPALM} &\multicolumn{2}{c|}{ALM-DQA-mod} \\
		\hline
		\multicolumn{1}{|c|}{Data} & \multicolumn{1}{|c|}{$m_0 | m_i$} & \multicolumn{1}{|c|}{$n_0 | n_i$}
		& \multicolumn{1}{|c||}{$N$}  & \multicolumn{1}{|c|}{Iter} &\multicolumn{1}{|c||}{Time(s)} &\multicolumn{1}{|c|}{Iter} & \multicolumn{1}{|c||}{Time (s)} &\multicolumn{1}{|c|}{Iter} & \multicolumn{1}{|c|}{Time (s)} \\ \hline
		
		\endhead
		
 qp-rand-m1-n20-N10-t1  &   1 $|$    1 &  20 $|$   20 &  10 & 321 &0.50 & 153 &8.43  &  12 &15.62 
 \\[3pt]\hline
 qp-rand-m50-n80-N10-t1  &  50 $|$   50 &  80 $|$   80 &  10 & 421 &0.64 & 268 &59.88  &  44 &192.17 
 \\[3pt]\hline
 qp-rand-m10-n20-N10-t2  &  10 $|$   10 &  20 $|$   20 &  10 &1501 &1.20 &2971 &208.56  &  54 &137.83 
 \\[3pt]\hline
 qp-rand-m50-n80-N10-t2  &  50 $|$   50 &  80 $|$   80 &  10 & 141 &0.20 &  92 &19.20  &  32 &150.32 
 \\[3pt]\hline
 tripart1  &2096 $|$  192 &2096 $|$ 2096 &  16 &1981 &3.01 &3880 &1212.86  &1422 &1113.26 
 \\[3pt]\hline
 tripart2  &8432 $|$  768 &8432 $|$ 8432 &  16 &6771 &51.65 &5000 &6369.20  &1610 &5723.31 
 \\[3pt]\hline
 qp-tripart1  &2096 $|$  192 &2096 $|$ 2096 &  16 & 653 &1.44 & 308 &94.60  & 114 &202.45 
 \\[3pt]\hline
 qp-tripart2  &8432 $|$  768 &8432 $|$ 8432 &  16 & 971 &9.79 & 347 &419.16  & 124 &1043.54 
 \\[3pt]\hline
 qp-pds1  &  87 $|$  126 & 372 $|$  372 &  11 & 971 &0.99 & 538 &49.18  & 535 &79.24 
 \\[3pt]\hline
 qp-SDC-r100-c50-l100-p1000-t1  &5000 $|$  150 &   0 $|$ 5000 & 100 &  32 &1.52 &  10 &37.49  &   7 &61.51 
 \\[3pt]\hline
 qp-SDC-r100-c50-l100-p1000-t2  &5000 $|$  150 &   0 $|$ 5000 & 100 &  31 &1.34 &   9 &34.22  &   2 &33.40 
 \\[3pt]\hline
 qp-SDC-r100-c50-l100-p5000-t1  &5000 $|$  150 &   0 $|$ 5000 & 100 &  32 &1.40 &  10 &37.75  &   8 &75.85 
 \\[3pt]\hline
 qp-SDC-r100-c50-l100-p5000-t2  &5000 $|$  150 &   0 $|$ 5000 & 100 &  31 &1.37 &   9 &34.50  &   3 &36.63 
 \\[3pt]\hline
 qp-SDC-r100-c50-l100-p10000-t1  &5000 $|$  150 &   0 $|$ 5000 & 100 &  32 &1.37 &  10 &37.93  &   9 &86.82 
 \\[3pt]\hline
 qp-SDC-r100-c50-l100-p10000-t2  &5000 $|$  150 &   0 $|$ 5000 & 100 &  31 &1.35 &   9 &34.78  &   3 &37.30 
 \\[3pt]\hline
 qp-SDC-r100-c100-l100-p1000-t1  &10000 $|$  200 &   0 $|$ 10000 & 100 &  31 &2.67 &  10 &73.50  &   7 &116.72 
 \\[3pt]\hline
 qp-SDC-r100-c100-l100-p1000-t2  &10000 $|$  200 &   0 $|$ 10000 & 100 &  31 &2.67 &   9 &68.08  &   2 &63.91 
 \\[3pt]\hline
 qp-SDC-r100-c100-l100-p5000-t1  &10000 $|$  200 &   0 $|$ 10000 & 100 &  31 &2.70 &  10 &74.02  &   8 &137.19 
 \\[3pt]\hline
 qp-SDC-r100-c100-l100-p5000-t2  &10000 $|$  200 &   0 $|$ 10000 & 100 &  31 &2.63 &   9 &68.04  &   2 &64.26 
 \\[3pt]\hline
 qp-SDC-r100-c100-l100-p10000-t1  &10000 $|$  200 &   0 $|$ 10000 & 100 &  32 &2.65 &  10 &74.65  &   8 &147.16 
 \\[3pt]\hline
 qp-SDC-r100-c100-l100-p10000-t2  &10000 $|$  200 &   0 $|$ 10000 & 100 &  31 &2.63 &   9 &67.66  &   3 &71.18 
 \\[3pt]\hline
 qp-SDC-r100-c100-l200-p20000-t1  &10000 $|$  200 &   0 $|$ 10000 & 200 &  41 &6.68 &  10 &183.21  &   8 &302.29 
 \\[3pt]\hline
 qp-SDC-r200-c100-l200-p20000-t1  &20000 $|$  300 &   0 $|$ 20000 & 200 &  34 &11.96 &  10 &360.48  &   7 &513.61 
 \\[3pt]\hline
 qp-SDC-r200-c200-l200-p20000-t1  &40000 $|$  400 &   0 $|$ 40000 & 200 &  31 &22.33 &  10 &783.49  &   7 &1068.04 
 \\[3pt]\hline
 M64-64  & 405 $|$   64 & 511 $|$  511 &  64 &1991 &3.16 &5000 &2874.88  & 625 &1241.12 
 \\[3pt]\hline

	\end{longtable}
\end{center}

%%%%%%%%%%%%%%%%%%%%%%%%%%%%%%%%%%%%
\section{A semi-proximal symmetric Gauss-Seidel based ADMM for the dual problem
(D)}

In the last section, we have designed the sPALM algorithm to solve
the primal problem (P)  directly.
One can also attempt to solve (P) via its dual problem (D) given in
\eqref{eq-struct-CCQP-dual}. Based on the structure in (D),
we find that it is highly conducive for us to employ a symmetric Gauss-Seidel based ADMM
(D) to
solve the problem, as we shall see later when the details are presented.

To derive the sGS-ADMM algorithm for solving
(D), it is more convenient for us to
 express (D) in a more compact form as follows:
\begin{eqnarray}
	\mbox{min}\;\{p(s)+f(y_{1:N},w,s)+q(z)+g(y_0,z) \mid \cF^*[y_{1:N};w;s] + \cG^*[y_0;z] = c\},
\label{eq-struct-CCQP-dual-compact}	
\end{eqnarray}
where $y_{1:N} = [y_1;\ldots; y_N]$, 	 and
\begin{eqnarray*}
\cF^* &:=&\big[\begin{array}{ccc} \cD^*, & -\bQ,& I  \end{array}\big], \quad
	\cG^*:=\big[\begin{array}{cc} \cA^*, & I  \end{array}\big],
	\\[5pt]
p(s) &:=& \mb{\theta}^*(-s),\;  f(y_{1:N},w,s):= -\inprod{b_{1:N}}{y_{1:N}} + \frac{1}{2}\inprod{w}{\bQ w}
+ \delta_{\mb{\cW}}(w) ,
\\[5pt]
q(z) &:=&\delta_{\bK}^*(-z),\; g(y_0,z):= -\inprod{b_0}{y_0}.
\end{eqnarray*}	
Here we take $\mb{\cW} = {\rm Range}(\bQ)$.
This is a multi-block linearly constrained convex
programming problem for which the direct application of the classical ADMM is
not guaranteed to converge. Thus we adapt the recently developed sGS-ADMM
\cite{CST,SCB-ADMM} whose convergence is guaranteed to solve
the dual problem (D).

Given a positive parameter $\sig$,
the augmented Lagrangian function for (D)
is given by
\begin{eqnarray*}
\begin{array}{lll}
 \cL_\sig(y,w,s,z;x)
 &=& p(s)+f(y_{1:N},w,s)+q(z)+g(y_0,z)  +
\\[5pt]
&&
\frac{\sig}{2}\norm{ \cF^*[y_{1:N};w;s]
+ \cG^*[y_0;z] -c + \frac{1}{\sig}x}^2 -\frac{1}{2\sig}\norm{x}^2
\nn\\[5pt]
&=&  \sum_{i=0}^N \theta_i^*(-s_i)+ \delta_{\cK_i}^*(-z_i)
+ \frac{1}{2}\inprod{w_i}{\cQ_i w_i} - \inprod{b_i}{y_i}
\nn \\[5pt]
&& + \frac{\sig}{2}
\norm{-\cQ_0 w_0 + \cA_0^* y_0 + s_0+z_0-c_0 +\sig^{-1}x_0}^2 -\frac{1}{2\sig}\norm{x_0}^2.
\nn \\[5pt]
&& + \sum_{i=1}^N \frac{\sig}{2}
\norm{-\cQ_i w_i + \cA_i^* y_0 +\cD_i^* y_i+ s_i+z_i-c_i +\sig^{-1}x_i}^2 -\frac{1}{2\sig}\norm{x_i}^2.
\end{array}
\end{eqnarray*}
Now to develop the sGS-ADMM, we need to analyze the block structure of the
quadratic terms in $\cL_\sig(y,w,s,z;x) $ corresponding the blocks $[y_{1:N};w;s]$ and $[y_0; z]$, which are respectively
given as follows:
%\begin{eqnarray*}
% \cF\cF^* &=& \left[\begin{array}{cc}
%  I &\bB^* \\[5pt] \bB & \bB\bB^*
% \end{array}\right]\;=\;
% \left[\begin{array}{cccc}
%  I & 0 & \cA_0^*  & 0  \\[5pt]
%  0 & I & \cA_{1:N}^* & \cD^* \\[5pt]
%  \cA_0 & \cA_{1:N} &\cA\cA^* &\cA_{1:N}\cD^* \\[5pt]
%  0 &\cD & \cD\cA_{1:N}^* & \cD\cD^*
% \end{array}\right]
% \\[5pt]
% &=&
% \underbrace{\left[\begin{array}{cccc}
%  0 & 0 & \cA_0^*  & 0  \\[5pt]
%  0 & 0& \cA_{1:N}^* & \cD^* \\[5pt]
%  0 & 0 &0  &\cA_{1:N}\cD^* \\[5pt]
%  0 &0& 0 & 0
% \end{array}\right]}_{\bU_\cF}
% +
% \underbrace{\left[\begin{array}{cccc}
%  I & 0 & 0  & 0  \\[5pt]
%  0 & I & 0 &0 \\[5pt]
%  0 & 0 &\cA\cA^*  &0 \\[5pt]
%  0 &0& 0 & \cD\cD^*
% \end{array}\right]}_{\bD_\cF} +\;\bU_\cF^*
% \\[5pt]
%\cG\cG^* & =& \left[\begin{array}{cc}
% I & -\bQ^* \\[5pt]
% -\bQ &\bQ\bQ^*
%\end{array} \right]  = \underbrace{\left[\begin{array}{cc}
% 0 & -\bQ^* \\[5pt]
% 0 & 0
%\end{array} \right] }_{\bU_\cG}+
%\underbrace{\left[\begin{array}{cc}
% I & 0  \\[5pt]
% 0 & \bQ\bQ^*
%\end{array} \right]}_{\bD_\cG} +\; \bU_\cG^*.
%\end{eqnarray*}
\begin{eqnarray*}
		\cF\cF^* &=& \left[\begin{array}{ccc}
			 \cD\cD^* & -\cD\bQ &\cD \\[5pt]
			 -\bQ\cD^* & \bQ^2 & -\bQ\\[5pt]
			     \cD^* & -\bQ & I
		\end{array}\right]
		\\[5pt]
		&=&
		\underbrace{\left[\begin{array}{ccc}
				0  & -\cD\bQ & \cD \\[5pt]
				0 & 0 & -\bQ \\[5pt]
				0 & 0 & 0
			\end{array}\right]}_{\bU_\cF}
		+
		\underbrace{\left[\begin{array}{ccc}
				\cD\cD^* & 0 & 0   \\[5pt]
				0 & \bQ^2 & 0 \\[5pt]
				0 & 0 & I
			\end{array}\right]}_{\bD_\cF} +\;\bU_\cF^*
		\\[5pt]
		\cG\cG^* & =& \left[\begin{array}{cc}
			\cA\cA^* & \cA\\[5pt]
			\cA^* & I
		\end{array} \right]  = \underbrace{\left[\begin{array}{cc}
			0 & \cA\\[5pt]
			0 & 0
		\end{array} \right] }_{\bU_\cG}+
	\underbrace{\left[\begin{array}{cc}
			\cA\cA^*& 0  \\[5pt]
			0 & I
	\end{array} \right]}_{\bD_\cG} +\; \bU_\cG^*.
\end{eqnarray*}
Based on the above (symmetric Gauss-Seidel) decompositions,
we define the following positive semidefinite linear operators
associated with the decompositions:
\begin{eqnarray}
 {\rm sGS}(\cF\cF^*) \;=\; \bU_\cF^*\bD_\cF^{-1}\bU_\cF, \quad
 {\rm sGS}(\cG\cG^*) \;=\; \bU_\cG^*\bD_\cG^{-1}\bU_\cG.
\end{eqnarray}
Note that here we view $\bQ$ as a linear operator defined on $\bW$ and because we take
$\bW = {\rm Range}(\bQ)$, $\bQ^2$ is positive definite on $\bW$
and hence $\bD_{\cF}$ is invertible. Since $\cA$ is assumed to have full row-rank,
$\bD_{\cG}$ is also invertible.

Given the current iterate $(y^k,s^k,w^k,z^k,x^k)$, the
basic template of the sGS-ADMM for \eqref{eq-struct-CCQP-dual-compact} at the $k$-th iteration
is given as follows.
\begin{description}
\item[Step 1.] Compute
\begin{eqnarray*}
 (y^{k+1}_{1:N},w^{k+1},s^{k+1}) \;=\;
 \mbox{argmin}_{y_{1:N},w,s}
 \left\{
 \begin{array}{l}
 p(s)+f(y_{1:N},w,s)
 \\[5pt]
 + \frac{\sig}{2}\norm{ \cF^*[y_{1:N};w;s] + \cG^*[y_0^k; z^k] -c + \frac{1}{\sig}x^k}^2
 \\[5pt]
 + \frac{\sig}{2}\norm{[y_{1:N};w;s]-[y_{1:N}^k;w^k; s^k]}^2_{{\rm sGS}(\cF\cF^*)}
 \end{array}
 \right\}.
\end{eqnarray*}
\item[Step 2.] Compute
\begin{eqnarray*}
 (y_0^{k+1},z^{k+1}) & = &\mbox{argmin}_{y_0,z}\left\{
 \begin{array}{l}
 q(z)+g(y_0,z) + \frac{\sig}{2}\norm{ \cF^*[y_{1:N}^{k+1};w^{k+1};s^{k+1}] + \cG^*[y_0;z] -c +\frac{1}{\sig} x^k}^2
 \\[5pt]
 + \frac{\sig}{2}\norm{[y_0;z]-[y_0^k;z^k]}^2_{{\rm sGS}(\cG\cG^*)}
 \end{array}
 \right\}.
\end{eqnarray*}
\item[Step 3.] Compute $x^{k+1} = x^k + \tau \sig (\cF^*[y_{1:N}^{k+1};w^{k+1};s^{k+1}] + \cG^*[y_0^{k+1}; z^{k+1}] - c)$,
where $\tau\in (0,\frac{1+\sqrt{5}}{2})$ is the steplength.
\end{description}

\medskip
By using the sGS-decomposition theorem in \cite{LSTb}, we can show that the
computation in Step 1 can be done by updating the blocks $(y_{1:N},w,s)$ in
a symmetric Gauss-Seidel fashion. Similarly, the computation in Step 2
can be done by updating the blocks $(y_0,z)$ in a symmetric Gauss-Seidel fashion.
With the above preparations, we can now give the detailed description of the
sGS-ADMM algorithm for solving \eqref{eq-struct-CCQP-dual}.

\medskip
\begin{description}
\item[sGS-ADMM on \eqref{eq-struct-CCQP-dual}.]
Given $(y^0,w^0,s^0,z^0,x^0)$
$\in \bY\times \bW \times \bX \times \bX \times \bX$,
 perform the following steps in
each iteration. Note that for notational convenience, we define $\cD_0=0$ in the algorithm.

\item[Step 1a.] Let $g^k = \cA^* y_0^k + z^k -c +\sig^{-1}x^k$.
Compute
\begin{eqnarray*}
 (\bar{y}_1^k,\ldots,\bar{y}_N^k) &=& \mbox{argmin}_{y_1,\ldots,y_N} \Big\{
\cL_\sig\big( (y_0^k,y_1,\ldots,y_N),w^k,s^k,z^k;x^k\big)
 \Big\},
\end{eqnarray*}
which can be done in parallel by computing for $i=1,\ldots,N$,
\begin{eqnarray*}
\bar{y}_i^{k} &=& \mbox{argmin}_{y_i} \Big\{ -\inprod{b_i}{y_i} + \frac{\sig}{2}\norm{ -\cQ_i w_i^k +\cD_i^* y_i + s^k_i + g_i^k}^2
\Big\}.
\end{eqnarray*}
Specifically, for $i=1,\ldots,N$, $\bar{y}_i^k$ is the solution of the following linear system:
\begin{eqnarray}
 \cD_i\cD_i^*\, y_i & = & \sig^{-1} b_i -\cD_i (-\cQ_i w_i^k + s^{k}_i + g^k_i).
  \label{eq-DDt}
\end{eqnarray}

\item[Step 1b]
Compute
$
 \bar{w}^k =\mbox{argmin}\big\{
 \cL_\sig\big( (y_0^k,\bar{y}_1^k,\ldots,\bar{y}_N^k),w,s^{k},z^k;x^k\big)\big\}
$
by computing in parallel for $i=0,1,\ldots,N,$
\begin{eqnarray*}
 \bar{w}_i^k &=& \mbox{argmin}_{w_i} \Big\{
 \frac{1}{2}\inprod{w_i}{\cQ_iw_i}
+\frac{\sig}{2}\norm{-\cQ_i {w}_i + \cD_i^* \bar{y}_i^{k}+ s^k_i +g^k_i}^2 \mid w_i \in \mbox{Range}(\cQ_i)
\Big\}.
\end{eqnarray*}
It is important to note that $\bar{w}^k_i$ is only needed theoretically
but not needed explicitly in practice. This is because in practical computation, only $\cQ_i \bar{w}^k_i$  is needed.
To compute $\cQ_i\bar{w}^k_i$, we first compute the solution $\tilde{w}_i^{k}$ of the
linear system below:
\begin{eqnarray}
(I + \sig\cQ_i )\tilde{w}_i \;=\; \sig (\cD_i^*\bar{y}_i^k +s_i^{k} + g^k_i).
\label{eq-Q}
\end{eqnarray}
Then
we can compute $\cQ_i \bar{w}^k_i = \cQ_i \tilde{w}_i^k$. The precise mechanism as to why the latter
equality is valid will be given in the remark after the presentation of this algorithm.

\item[Step 1c.] Compute
\begin{eqnarray*}
 (s_0^{k+1},\ldots,s_N^{k+1}) &=& \mbox{argmin}_{s_0,\ldots,s_N} \Big\{
\cL_\sig( (y_0^k,\bar{y}_1^{k},\ldots,\bar{y}_N^{k}),\bar{w}^k,(s_0,s_1,\ldots,s_N),z^k;x^k)
 \Big\},
\end{eqnarray*}
which can be done in parallel by computing for $i=0,1,\ldots,N,$
\begin{eqnarray*}
s_i^{k+1} &=&
 \mbox{argmin}_{y_i} \Big\{ \theta_i^*(-s_i) + \frac{\sig}{2}\norm{
 -\cQ_i\bar{w}_i^{k} +\cD_i^* \bar{y}^{k}_i + s_i +g^k_i}^2
\Big\}
\\[5pt]
&=&- {\rm Prox}_{\theta_i^*/\sig}(-\cQ_i \bar{w}_i^{k} + \cD_i^* \bar{y}_i^{k}
 + g_i^k)
\\[5pt]
&=& \frac{1}{\sig} \mbox{Prox}_{\sig \theta_i}\big(\sig(-\cQ_i \bar{w}_i^k + \cD_i^* \bar{y}_i^k + g_i^k)\big)
-(-\cQ_i\bar{w}_i^k+ \cD_i^* \bar{y}_i^k + g_i^k).
\end{eqnarray*}

\item[Step 1d] Compute
$
 w^{k+1} =\mbox{argmin}\big\{
 \cL_\sig\big( (y_0^k, \bar{y}_1^k,\ldots,\bar{y}_N^k),w,s^{k+1},z^{k};x^k\big)  \big\}
$
by computing in parallel for $i=0,1,\ldots,N,$
\begin{eqnarray*}
 {w}_i^{k+1} &=& \mbox{argmin}_{w_i} \Big\{
 \frac{1}{2}\inprod{w_i}{\cQ_iw_i}
+\frac{\sig}{2}\norm{-\cQ_i {w}_i+ \cD_i^* \bar{y}_i^k + s_i^{k+1} +g^k_i}^2 \mid w_i \in \mbox{Range}(\cQ_i)
\Big\}.
\end{eqnarray*}
Note that the same remark in Step 1b is applicable here.

\item[Step 1e] Compute
\begin{eqnarray*}
 (y_1^{k+1},\ldots,y_N^{k+1} ) &=& \mbox{argmin}_{y_1,\ldots,y_N} \Big\{
\cL_\sig\big( (y_0^k,y_1,\ldots,y_N),w^{k+1},s^{k+1},z^k;x^k\big)
 \Big\},
\end{eqnarray*}
which can be done in parallel by computing for $i=1,\ldots,N$,
\begin{eqnarray*}
 {y}_i^{k+1} &=& \mbox{argmin}_{y_i} \Big\{ -\inprod{b_i}{y_i} + \frac{\sig}{2}\norm{ -\cQ_i w_i^{k+1} +\cD_i^* y_i + s^{k+1}_i + g^k_i}^2
\Big\}.
\end{eqnarray*}

%%
%% Step 2
%%
\item[Step 2a.] Let $\h^k = -\bQ w^{k+1} +\cD^* y^{k+1} + s^{k+1} - c +\sig^{-1} x^k.$
Compute
\begin{eqnarray*}
 \bar{y}_0^{k} &=& \mbox{argmin}_{y_0} \Big\{
\cL_\sig\big( (y_0,y_1^{k+1},\ldots,y_N^{k+1}),w^{k+1},s^{k+1},z^k;x^k\big)
 \Big\}
\\[5pt]
&=& \mbox{argmin}_{y_0} \Big\{
-\inprod{b_0}{y_0} +
\frac{\sig}{2}\norm{ \cA_0^* y_0 + z_0^k + h^k_0}^2
+ \sum_{i=1}^N \frac{\sig}{2}\norm{ \cA_i^* y_0+ z^k_i +h^k_i}^2
\Big\}.
\end{eqnarray*}
Specifically, $\bar{y}_0^k$ is the solution to the following linear system of equations:
\begin{eqnarray}
 \Big(\sum_{i=0}^N \cA_i\cA_i^* \Big) y_0 & = &
 \sig^{-1} b_0 - \sum_{i=0}^N \cA_i ( z_i^k + h_i^k).
 \label{eq-AAt}
\end{eqnarray}

\item[Step 2b] Compute $
 z^{k+1} = \mbox{argmin}\big\{
 \cL_\sig\big( (\bar{y}_0^k, y_1^{k+1},\ldots,y_N^{k+1}), w^{k+1},s^{k+1},z;x^k\big)\big\}$
 by computing in parallel for $ i=0,1,\ldots,N$,
\begin{eqnarray*}
 z_i^{k+1} &=& \mbox{argmin}_{z_i} \Big\{
  \delta_{\cK_i}^*(-z_i)
+\frac{\sig}{2}\norm{\cA_i^* \bar{y}_0^k +  z_i  + h^k_i}^2
\Big\}
\;=\; -\mbox{Prox}_{\sig^{-1}\delta_{\cK_i}^*}(\cA_i^*\bar{y}_0^k + h_i^k)
\\[5pt]
&=& \frac{1}{\sig} \Pi_{\cK_i}\big(\sig (\cA^*_i \bar{y}_0^k + h_i^k) \big)
- (\cA^*_i \bar{y}_0^k + h_i^k).
\end{eqnarray*}

\item[Step 2c]
Compute
\begin{eqnarray*}
&& \hspace{-0.7cm}
 {y}_0^{k+1} \;=\; \mbox{argmin}_{y_0} \Big\{
\cL_\sig\big( (y_0,y_1^{k+1},\ldots,y_N^{k+1}),w^{k+1},s^{k+1},z^{k+1};x^k\big)
 \Big\}
\\[5pt]
&=&
 \mbox{argmin}_{y_0} \Big\{
-\inprod{b_0}{y_0} +
\frac{\sig}{2}\norm{ \cA_0^* y_0 + z^{k+1}_0 + h_0^k}^2
+ \sum_{i=1}^N \frac{\sig}{2}\norm{ \cA_i^* y_0 + z^{k+1}_i + h_i^k}^2
\Big\}.
\end{eqnarray*}
Note that the computation in Step 2a is applicable here.

\item[Step 3] Compute
\begin{eqnarray*}
x^{k+1} &=& x^k + \tau \sigma (-\bQ w^{k+1} + \bB^* y^{k+1}+ s^{k+1} + z^{k+1} - c),
\end{eqnarray*}
where $\tau\in (0,\frac{1+\sqrt{5}}{2})$ is the steplength.
\end{description}

\bigskip
Now we make some important remarks concerning the computations in sGS-ADMM.
\begin{enumerate}
\item If the term $\mb{\theta}\equiv 0$ in Step 1c, then this step is vacuous, and
Step 1b and Step 1d are identical. Hence the computation needs only to be
done for Step 1d. Hence Step 1 only consists of Step 1a, 1d, and 1e.

\item If $\mb{\cQ} \equiv 0$, then Step 1b and 1d are vacuous. Therefore
Step 1 only consists of Step 1a, 1c, and 1e.

\item The computation in Step 1d can be
omitted if the quantity $\bar{w}_i^k$ computed in Step 1b
is already a sufficiently good approximate solution to the current subproblem.
More precisely, if the approximation $\bar{w}_i^k$ for
$w_i^{k+1}$  satisfies the admissible accuracy condition
required in the inexact sGS-ADMM designed in \cite{CST}, then we
can just set $w_i^{k+1} = \bar{w}_i^k$ instead of using the
exact solution to the current subproblem.
Similar remark is also applicable to the computation in Step 1e and Step 2c.

\item The sGS-ADMM in fact has the flexibility of allowing for
inexact computations as already shown in \cite{CST}.
 While the computation in Step 1a and 1e (similarly for Step 1b and 1d, Step 2a and 2c)
 are assumed to be done exactly
(up to machine precision), the computation can in fact be done inexactly
subject to a certain predefined accuracy requirement on the
computed approximate solution.
Thus iterative methods such as the preconditioned conjugate gradient (PCG) method can be
used to solve the linear systems when their dimensions are too large.
We omit the details here for
the sake of brevity.

\item
In solving the linear system \eqref{eq-AAt}, the $m_0\times m_0$ symmetric positive
definite matrix $\sum_{i=0}^N \cA_i\cA_i^*$ is fixed, and one can pre-compute
the matrix if it can be stored in the memory and its Cholesky factorization can be
computed at a reasonable cost. Then in each sGS-ADMM iteration, $\bar{y}_0^k$ and $y_0^{k+1}$ can
be computed cheaply by solving triangular linear systems.
In the event when computing the coefficient matrix or its Cholesky factorization is out of reach,
one can use a PCG method to solve the linear system. In
that case, one can implement the computation of the matrix-vector product in parallel by
computing $\cA_i \cA_i^* y_0$ in parallel for $i=0,1,\ldots,N$, given any $y_0$.
Note that when the PCG method is employed, the use of the  inexact sGS-ADMM framework just mentioned above
will become necessary.

The same remark above
also applies to the linear system \eqref{eq-DDt} for each $i=1,\ldots,N$.

For the multi-commodity flow problem which we will consider
later in the numerical experiments, we note that  the linear system in
\eqref{eq-AAt} has a very simple coefficient matrix given by $\sum_{i=0}^N \cA_i\cA_i^* = (N+1) I_m$, and the coefficient matrix $\cD_i\cD_i^*$ in \eqref{eq-DDt} is
equal to the Laplacian matrix of the network graph for all $i=1,\ldots,N.$
Thus both \eqref{eq-AAt} and \eqref{eq-DDt} can be solved efficiently
by a direct solver.

\item
In Step 1b, we claimed that $\cQ_i \bar{w}^k_i = \cQ_i \tilde{w}^k_i$. Here we show why the result
holds. For simplicity, we assume that $\cQ_i$ is a symmetric positive semidefinite matrix rather than a
linear operator. Consider the spectral decomposition $\cQ_i = UDU^T$, where $D\in \Re^{r\times r}$ is a
diagonal matrix whose diagonal elements are the positive eigenvalues of $\cQ_i$ and
the columns of
$U\in \Re^{m_i\times r}$ are their corresponding orthonormal set of eigenvectors.
We let $V\in \Re^{m_i\times (m_i -r)}$ be the matrix whose columns  form
an orthonormal set of eigvectors of $\cQ_i$ correspond to the zero eigenvalues.
%Based on
With this decomposition and the parameterization $w_i = U \xi$ (because $w_i \in {\rm Range}(\cQ_i)$),
the minimization for $\bar{w}^k_i$ is equivalent to the following:
\begin{eqnarray}
\mbox{argmin} \Big\{\frac{1}{2}\inprod{\xi}{D\xi} + \frac{\sig}{2}\norm{ D \xi - U^T g}^2 + \frac{\sig}{2}
\norm{V^Tg}^2 \mid \xi \in \Re^r
\Big\},
\label{eq-D}
\end{eqnarray}
where we have set $g = z_i^k+h_i^k$ for convenience.
Now from solving \eqref{eq-Q}, we get that
\begin{eqnarray*}
 (I + \sig D) U^T \tilde{w}_i^k = \sig U^T g, \quad V^T \tilde{w}^k_i = \sig V^T g.
\end{eqnarray*}
This show that $U^T\tilde{w}^k_i$ is the unique solution to the problem \eqref{eq-D}.
Hence $\bar{w}^k_i = U (U^T \tilde{w}^k_i)$ is the unique solution to \eqref{eq-Q}.
From here, we have that $\cQ_i \bar{w}^k_i = UDU^T (UU^T \tilde{w}^k_i) = UDU^T \tilde{w}^k_i =  \cQ_i \tilde{w}^k_i$.

\end{enumerate}

%%%%%%%%%%%%%%%%%%%%%

\subsection{Convergence theorems of sGS-ADMM}

The convergence theorem of sGS-ADMM can be established directly by using known results from \cite{CST} and \cite{ZWZ}. Here we present  the global convergence result and the linear rate of convergence for the convenience of reader.

In order to state the convergence theorems, we need some definitions.
\begin{defn}
	Let $\mathbb{F}:\cX \rightrightarrows\cY$ be a multivalued mapping and denote its inverse by $\mathbb{F}^{-1}$. The graph of multivalued function $\mathbb{F}$ is defined by
	${\rm gph}\mathbb{F}:=\{(x,y)\in\cX\times \cY\mid y\in \mathbb{F}(x)\}$.
\end{defn}

Denote $u:=(y,w,s,z,x)\in\mathcal{U}:=\bY\times\bW\times \bX\times \bX\times \bX$. The KKT mapping $\mathcal{R}:\mathcal{U}\rightarrow\mathcal{U}$ of \eqref{eq-CCQP} is defined by
\begin{eqnarray}
\mathcal{R}(u):=\begin{pmatrix}
\bB x - b \\
-\bQ w + \bB^* y + s + z - c \\
\bQ w - \bQ x \\
x - {\rm Prox}_{\mb{\theta}} (x-s) \\
x - \Pi_{\bK}(x-z) \\
\end{pmatrix}.
\label{eq-KKTmapping}
\end{eqnarray}

Denote the set of KKT points by $\bar{\Omega}$. The  KKT mapping $\mathcal{R}$ is said to be
metrically subregular at $(\bar{u},0)\in {\rm gph}\mathcal{R}$ with modulus $\eta>0$ if
 there exists a scalar $\rho>0$ such that
$${\rm dist}(u,\bar{\Omega})\le \eta\|\mathcal{R}(u)\| \quad \forall u \in \{u\in \mathcal{U}: \|u-\bar{u}\|\le \rho\}.$$

Now we are ready to present the convergence theorem of sGS-ADMM.
\begin{theorem}
	Let $\{u^k:=(y^k,w^k,s^k,z^k;x^k)\}$ be the sequence generated by sGS-ADMM. Then, we have the following results.
	\\[5pt]
	(a) The sequence $\{(y^k,w^k,s^k,z^k)\}$ converges to an optimal solution of the compact
	form \eqref{eq-CCQP-dual} of the dual problem (D),
	and the sequence $\{x^k\}$ converges to an optimal solution of the compact form \eqref{eq-CCQP}
	of the primal problem (P).
	\\[5pt]
	(b) Suppose that the sequence $\{u^k\}$ converges to a KKT point $\bar{u}:=(\bar{y}^k,\bar{w}^k,\bar{s}^k,\bar{z}^k,\bar{x}^k)$ and the KKT mapping $\mathcal{R}$ is metrically subregular at $(\bar{u},0)\in {\rm gph}\mathcal{R}$. Then the sequence $\{u^k\}$ is linearly convergent to $\bar{u}$.
\end{theorem}

\begin{Proof}
	(a) The global convergence result follows from that in \cite{CST}.
	(b)
	The result follows directly by applying the convergence result in \cite[Proposition 4.1]{ZWZ} (which slightly improves
	an earlier result in \cite{HSZ})  to the compact formulation \eqref{eq-CCQP-dual} of (D).
\end{Proof}

\begin{remark} By Theorem 1 and Remark 1 in \cite{LSTc}, we know that when (P) is
a convex  programming problem where for each $i=0,\ldots,N$,  $\theta_i$ is piecewise linear-quadratic or strongly convex, and
$\cK_i$ is polyhedral, then
$\cal R$ is metrically
subregular at $(\bar{u},0)\in {\rm gph}\mathcal{R}$  for any KKT point $\bar{u}.$
Thus sGS-ADMM converges locally at a linear rate to an optimal solution of (P) and (D)
under the previous conditions on $\theta_i$ and $\cK_i$. In particular,
for the special case of a primal block angular quadratic programming problem where $\theta_i \equiv  0$ and $\cK_i  = \mathbb{R}^{n_i}_+$ for all $i$, we know that sGS-ADMM is locally linearly convergent, which can even be proven  to converge   globally  linearly.
\end{remark}

%%%%%%%%%%%%%%%%%%%%%
\subsection{Computational cost}
Now we would discuss the main computational cost of sGS-ADMM. We could observe that the most time-consuming computations are in solving large linear system of equations
in Step 1a, 1b, 1d, 1e, 2a, and 2c.

In general, suppose for every iteration we need to solve a $d\times d$ linear system of equations:
\begin{eqnarray}
Mx = r.
\end{eqnarray}
Assuming that $M$ is stored, then we can compute its Cholesky factorization at the cost of $O(d^3)$ operations, which needs only to be done once at the very beginning of the algorithm. After that, whenever we need to solve the equation, we just need to compute the right-hand-side vector $r$ and solve two $d\times d$ triangular systems of linear equations at the cost of $O(d^2)$ operations.

We can roughly summarize the costs incurred in solving $Mx=r$ as follows:
\begin{itemize}
\item[$(C_1)$] Cost for computing the coefficient matrix $M$ (only once at the beginning of algorithm);

\item[$(C_2)$] Cost for computing Cholesky factorization of $M$ (only once at the beginning of algorithm);

\item[$(C_3)$] Cost for computing right-hand-side vector $r$;

\item[$(C_4)$] Cost for solving two triangular systems of linear equations.
\end{itemize}

The computational cost $C_1,C_2,C_3,C_4$ above for each of the equations in Step 1a, 1b, 1d, 1e, 2a, and 2c are tabulated in Table \ref{table-cost}.

\begin{center}
	\setlength{\tabcolsep}{6pt}
	\begin{longtable}{| c || c | c | c | c |}
			\caption{\footnotesize
			Computational cost for solving the linear systems of equations in each of the steps. }
		\label{table-cost}
		\\
		\hline
		\multicolumn{1}{|c||}{Step} & \multicolumn{1}{|c|}{$C_1$ (once)} & \multicolumn{1}{|c|}{$C_2$ (once)} & \multicolumn{1}{|c|}{$C_3$ (each iteration)} & \multicolumn{1}{|c|}{$C_4$ (each iteration)} \\[3pt]
		\hline
		$\begin{array}{c}
		\mbox{1a and 1e} \\
		(i=1,\ldots,N)
		\end{array}$
		& $O(m_i^2 n_i)$ & $O(m_i^3)$ & $O(n_i^2+m_i n_i)$ & $O(m_i^2)$ \\[3pt] \hline
				$\begin{array}{c}
		\mbox{1b and 1d} \\
		(i=1,\ldots,N)
		\end{array}$
		& $O(n_i^2)$ & $O(n_i^3)$ & $O(m_i n_i)$ & $O(n_i^2)$ \\[3pt] \hline
		2a and 2c & $O(m_0^2 n_0)$ & $O(m_0^3)$ & $O(m_0 n_0)$ & $O(m_0^2)$ \\[3pt] \hline
		
	\end{longtable}
\end{center}

\section{Numerical experiments}
\label{sec-5}

In this section, we evaluate the performance of the algorithm we have designed for solving
the problem (P). We conduct numerical experiments on three major types of
primal block angular model, including linear, quadratic, and nonlinear problems. Apart from randomly generated datasets, we would demonstrate that our algorithms can be quite efficient in solving realistic problems encountered in the literature.

\subsection{Stopping condition}
Based on the optimality conditions in \eqref{KKT}, we measure the accuracy of a computed solution by the following relative residuals:
\[
\eta = \max\{\eta_P,\eta_D,\eta_Q,\eta_{K},\eta_{S}\},
\]
where
\begin{eqnarray*}
	\eta_P &=& \frac{\norm{\bB x-b}}{1+\norm{b}}, \quad
	\eta_D = \frac{\norm{-\bQ w + \bB^*y + s+z - c}}{1+\norm{c}}, \quad
	\eta_Q = \frac{\norm{\bQ w - \bQ x}}{1+\norm{\bQ}}, \\
	\eta_{K} &=& \frac{\norm{x-\Pi_{\bK}(x-z)}}{1+\norm{x}+\norm{z}}, \quad
	\eta_{S} = \frac{\norm{x-{\rm Prox}_\theta (x-s)}}{1+\norm{x}+\norm{s}}.
	%\eta_{K_1} &=& \frac{\norm{x-{\Pi_{\bK}}(x)}}{1+\norm{x}}, \quad
	%\eta_{K_3} = \frac{\|z-{\Pi_{\bK^*}}(z)\|}{1+\norm{z}},
\end{eqnarray*}
%In addition, we also compute the duality gap by:
%\[\eta_{gap} = \frac{|\text{obj}_\text{Primal} - \text{obj}_\text{Dual}|}{1+|\text{obj}_\text{Primal}|+|\text{obj}_\text{Dual}|},\]
%where $\text{obj}_\text{Primal}:=\mb{\theta}(x) + \frac{1}{2}\inprod{x}{\bQ x} + \inprod{c}{x}$ and $\text{obj}_\text{Dual}:=  \inprod{b}{y}-\frac{1}{2}\inprod{w}{\bQ w} -\mb{\theta^*}(-s)  -\delta_{\bK}^*(-z)$
%are the primal and dual objective functions respectively.
We terminate our algorithm when $\eta \le 10^{-5}$.

%%%%%%%%%%%%%%%%%%%%%%%%%%%%%%%%%%%
\subsection{Block angular  problems with linear objective functions}
\label{subsec-linear}

In this subsection, we perform numerical experiments on minimization
problems having linear objective functions and primal block angular constraints. Multicommodity flow (MCF) problems are one of the main representative in this  class of problems. It is a model to solve the routing problem of multiple commodities throughout a network from a set of supply nodes to a set of demand nodes. These problems usually exhibit primal block angular structures due to the network nature in the constraints.

Consider a connected network graph $(\cN,\cE)$ with $m$ nodes and $n=|\cE|$ arcs
for which $N$ commodities must be transported through the network. We assume that
each commodity has a single source-sink pair $(s_k,t_k)$ and we are given
the flow $r_k$ that must be transported from $s_k$ to $t_k$, for $k=1,\ldots,N$.
Let $M\in \mathbb{R}^{m\times |\cE|}$ be node-arc incidence matrix of the graph.
Then the MCF problem can be expressed in the form given in (P) with the following data:
\begin{eqnarray*}
     && \cK_0 = \{x_0 \in \mathbb{R}^{n} \mid 0 \leq x_0 \leq u\}, \quad
      \cK_i = \mathbb{R}^{n}_+, \;\; i=1,\ldots, N , \\
	&&\cQ_i = 0, \;\; \theta_i(\cdot) = 0, \;\; \forall i = 0,1,...,N,\\
	&& \cA_0 = I_n, \;\; \cA_i = -I_n, \;\; \forall i = 1,...,N,\\
	&&\cD_1 = \cD_2 = \cdots = \cD_N = M \text{ is the node-arc incidence matrix.}
\end{eqnarray*}
For this problem, $x_i$ denotes the flow of the $i$-th commodity ($i=1,\ldots,N$) through the
network,  $x_0$ is the total flow, and $u$ is a given upper bound vector on the total flow.

\subsubsection{Description of datasets}
Following \cite{Castro2011}, the datasets we used are as follows.
\begin{description}
	\item[\textbf{tripart} and \textbf{gridgen}:] These are five multicommodity instances obtained with the Tripart and Gridgen generators. They could be downloaded from \\ \url{http://www-eio.upc.es/~jcastro/mmcnf_data.html}.
	
	\item[\textbf{pds}:] The PDS problems come from a model of transporting patients away from a place of military conflict. It could be downloaded from \\ \url{http://www.di.unipi.it/optimize/Data/MMCF.html#Pds}.
	
	\item[\textbf{M\{$n$\}-\{$k$\}}:] These are the problems generated by the Mnetgen generator, which is one of the most famous random generator of Multicommodity Min Cost Flow instances. Here $n$ is the number of nodes in the network and $k$ is the number of commodity. It could be downloaded from \url{http://www.di.unipi.it/optimize/Data/MMCF.html#MNetGen}.
\end{description}

%%%%%%%%%%%%%%%%%%%%%%%
\subsubsection{Numerical results}
In Table \ref{table-linear}, we compare our sGS-ADMM algorithm against the solvers Gurobi and BlockIP.
We should emphasize that Gurobi is a state-of-the-art solver for solving general linear and quadratic programming problems. Although it is not a specialized algorithm for primal block angular problems, it has been so powerful in solving sparse general linear and convex quadratic programming problems that it should be used as the benchmark for any newly developed algorithm.
On the other hand, BlockIP \cite{Castro2016} is an efficient interior-point algorithm specially designed for solving primal block angular problems, especially those arising from MCF problems. As reported in \cite{Castro2016}, it has been successful in solving many large scale instances of primal block angular LP and QP problems.

In the following numerical experiments, we employ Gurobi directly on the compact formulation \eqref{eq-CCQP}. To be more specific,
we input $\cB$ as a general sparse matrix.
The feasibility and objective gap tolerance is set to be 1e-5, and the number of threads is set to be 1. All the other parameters remain as default setting. Similarly for BlockIP, all the three tolerances (primal and dual feasibility, and relative objective gap) are set to be 1e-5 for consistency. Its maximum number of iteration is set to be 500.

\begin{center}
\begin{scriptsize}
	\setlength{\tabcolsep}{4pt}
	\begin{longtable}{| l | c | c | c || c | c || c | c || c | c |}	
			\caption{\footnotesize
			Comparison of computational results between sGS-ADMM, Gurobi, and BlockIP for \textbf{linear} primal block angular problems. All the results are obtained using a \textbf{single thread}. `Iter' under the column for Gurobi means
			the total number of simplex iterations.}
		\label{table-linear}
		\\
		\hline
		\multicolumn{4}{|c||}{} & \multicolumn{2}{c||}{sGS-ADMM} & \multicolumn{2}{c||}{Gurobi} & \multicolumn{2}{c|}{BlockIP} \\
		\hline
		\multicolumn{1}{|c|}{Data} & \multicolumn{1}{|c|}{$m_0 | m_i$} & \multicolumn{1}{|c|}{$n_0 | n_i$}
		& \multicolumn{1}{|c||}{$N$}  & \multicolumn{1}{|c|}{Iter} &\multicolumn{1}{|c||}{Time(s)} &\multicolumn{1}{|c|}{Iter} & \multicolumn{1}{|c||}{Time (s)} &\multicolumn{1}{|c|}{Iter} & \multicolumn{1}{|c|}{Time (s)} \\ \hline
		
		\endhead
		
 tripart1  &2096 $|$  192 &2096 $|$ 2096 &  16 &1981 &3.01 &5155 &0.78 &  48 &1.23
 \\[3pt]\hline
 tripart2  &8432 $|$  768 &8432 $|$ 8432 &  16 &6771 &51.65 &42070 &42.81 &  67 &10.32
 \\[3pt]\hline
 tripart3  &16380 $|$ 1200 &16380 $|$ 16380 &  20 &5561 &104.96 &85390 &189.37 &  81 &48.70
 \\[3pt]\hline
 tripart4  &24815 $|$ 1050 &24815 $|$ 24815 &  35 &8581 &343.32 &246340 &1685.50 & 115 &139.36
 \\[3pt]\hline
 gridgen1  &3072 $|$ 1025 &3072 $|$ 3072 & 320 &7541 &409.75 &497709 &8039.40 & 203 &1589.04
 \\[3pt]\hline
 pds15  &1812 $|$ 2125 &7756 $|$ 7756 &  11 &2893 &22.60 &8545 &1.01 &  81 &12.19
 \\[3pt]\hline
 pds30  &3491 $|$ 4223 &16148 $|$ 16148 &  11 &4471 &111.49 &27645 &4.79 & 110 &51.66
 \\[3pt]\hline
 pds60  &6778 $|$ 8423 &33388 $|$ 33388 &  11 &7719 &465.06 &70168 &17.57 & 145 &403.23
 \\[3pt]\hline
 pds90  &8777 $|$ 12186 &46161 $|$ 46161 &  11 &5315 &479.59 &100858 &25.18 & 162 &822.45
 \\[3pt]\hline
 M64-64  & 405 $|$   64 & 511 $|$  511 &  64 &1991 &3.16 &7601 &0.77 &  51 &0.85
 \\[3pt]\hline
 M128-64  & 936 $|$  128 &1171 $|$ 1171 &  64 &2601 &7.32 &18108 &3.93 &  52 &3.15
 \\[3pt]\hline
 M128-128  & 979 $|$  128 &1204 $|$ 1204 & 128 &3801 &28.89 &32736 &7.24 & 127 &11.75
 \\[3pt]\hline
 M256-256  &1802 $|$  256 &2204 $|$ 2204 & 256 &6821 &225.31 &103561 &18.53 &  97 &89.92
 \\[3pt]\hline
 M512-64  &3853 $|$  512 &4768 $|$ 4768 &  64 &2631 &45.77 &48235 &8.76 &  72 &48.44
 \\[3pt]\hline
 M512-128  &3882 $|$  512 &4786 $|$ 4786 & 128 &3581 &137.23 &87659 &17.96 &  97 &144.77
 \\[3pt]\hline
 M512-512  & 707 $|$  512 &1797 $|$ 1797 & 512 &7021 &373.58 &199260 &16.79 & 146 &308.95
 \\[3pt]\hline

	\end{longtable}
\end{scriptsize}	
\end{center}

From Table \ref{table-linear}, we observe that Gurobi is the fastest to solve 11 out of 16 instances.
Gurobi is extremely fast in solving the {\tt pdsxx} and {\tt Mxxx-xx} problems but have
difficulty in solving {\tt tripart4} and {\tt gridgen1} efficiently. On the other hand,
sGS-ADMM and BlockIP are highly efficient in solving the latter instances.
%Comparing sGS-ADMM and BlockIP, the former is more efficient in solving
%three of the largest instances.
On the other hand, BlockIP is the fastest when solving the {\tt tripart2,3,4} instances
while sGS-ADMM is the fastest in solving the {\tt gridgen1} and {\tt M512-128} instances.

Our sGS-ADMM solver outperforms Gurobi when the instance is both hard and huge, for example, {\tt tripart4} and {\tt gridgen1}.
For the latter instance, it is in fact the fastest solver.
We also noticed that BlockIP is quite sensitive to the practical setting of the upper bound on the unbounded variables. For example,
setting ``9e6'' and ``9e8'' as the upper bounds for  the unbounded variables can lead to a
significant difference in the number of iterations.

%%%%%%%%%%%%%%%%%%%%%
\subsection{Block angular problems with convex quadratic objective functions}

In this subsection, we perform numerical experiments on optimization problems having convex quadratic objective functions and primal block angular constraints.

One of the main class of this type of problem is again  from the multicommodity flow problem.  Following \cite{Castro2016}, we add in the quadratic objective term, $\cQ_i = 0.1I, \; \forall i = 0,...,N$.
The corresponding datasets start with a prefix \texttt{"qp-"}, including \texttt{tripart}, {\tt gridgen} and \texttt{pds}.

Another main class of quadratic primal block angular problems arises in the field of statistical disclosure control. Castro \cite{Castro2005} studied the controlled tabular adjustment (CTA) to find a closest, perturbed, yet safe table given a three-dimensional table for which the content need to be protected. In particular, we have
\begin{eqnarray*}
	&&\cQ_i = I, \; \theta_i(\cdot) = 0, \; i=0,\ldots,N,\\
&& 	\cA_0 = I, \; \cA_i = -I, \; \forall i = 1,...,N,\\
	&&\cD_1 = \cD_2 = \cdots = \cD_N \text{ is a node-arc incidence matrix}
\end{eqnarray*}
and $\cK_i$ ($i=0,1,\ldots,N$) is the same as in section \ref{subsec-linear}.

\subsubsection{Description of datasets}
The datasets we used are as follows.
\begin{description}
	\item[rand:] These instances are randomly generated sparse problems.
	Here we generated two types of problems.
	\begin{itemize}
		\item Type 1 problem (with suffix \texttt{-t1}) has diagonal quadratic objective cost, i.e. $\cQ_i$ is a random diagonal matrix given by {\tt spdiags(rand(n\_i,1),0,n\_i,n\_i)}.
		\item Type 2 problem (with suffix \texttt{-t2}) does not necessarily have diagonal quadratic objective cost. In this case $\cQ_i$ is still very sparse but remained to be positive semidefinite. We use the following routine to generate $\cQ_i$ for every $i=0,1,...,N$:
		{\center \tt tmp=sprandn(n\_i,n\_i,0.1); $Q_i$ = tmp*tmp'.}
	\end{itemize}

	For both types of problems, we generate $\cA_i$ and $\cD_i$ similarly for $i=0,...,N$ using {\sc Matlab} command \texttt{sprandn} with density 0.5 and 0.3 respectively. Note that by convention we have $\cD_0=0$.
	
	\item[L2CTA3D:] This is an extra large instance (with a total of 10M variables and 210K constraints) provided in  \url{http://www-eio.upc.es/~jcastro/huge_sdc_3D.html}.
	
	\item[SDC:] These are some of the CTA instances we generated using the generator
	provided by J. Castro at \url{http://www-eio.upc.es/~jcastro/CTA_3Dtables.html}.
\end{description}

\subsubsection{Numerical results}
As in the last subsection, we compare our sGS-ADMM algorithm against Gurobi and BlockIP solver in  Table \ref{table-quad}.
\begin{center}
	\setlength{\tabcolsep}{3pt}
	\scriptsize
	\begin{longtable}{| l| c | c | c || c | c || c | c || c | c |}	
			\caption{\footnotesize
			Comparison of computational results between sGS-ADMM, Gurobi, and BlockIP for \textbf{quadratic} primal block angular problems. All the results are obtained using \textbf{single thread}.  `Iter' under the column for Gurobi means the total number of  barrier iterations. A `/' under the column for BlockIP means that the solver runs out of memory, and a `*' means the solver is not compatible to solve the problem.}
		\label{table-quad}
		\\
		\hline
		\multicolumn{4}{|c||}{} & \multicolumn{2}{c||}{sGS-ADMM} & \multicolumn{2}{c||}{Gurobi} & \multicolumn{2}{c|}{BlockIP} \\
		\hline
		\multicolumn{1}{|c|}{Data} & \multicolumn{1}{|c|}{$m_0 | m_i$} & \multicolumn{1}{|c|}{$n_0 | n_i$}
		& \multicolumn{1}{|c||}{$N$}  & \multicolumn{1}{|c|}{Iter} &\multicolumn{1}{|c||}{Time(s)} &\multicolumn{1}{|c|}{Iter} & \multicolumn{1}{|c||}{Time (s)} &\multicolumn{1}{|c|}{Iter} & \multicolumn{1}{|c|}{Time (s)} \\ \hline
		
		\endhead
		
 qp-rand-m50-n80-N10-t1  &  50 $|$   50 &  80 $|$   80 &  10 & 421 &0.64 &  14 &0.28 &  29 &0.13
 \\[3pt]\hline
 qp-rand-m1000-n1500-N10-t1  &1000 $|$ 1000 &1500 $|$ 1500 &  10 & 748 &57.67 &  15 &1641.91 &  39 &360.12
 \\[3pt]\hline
 qp-rand-m100-n200-N100-t1  & 100 $|$  100 & 200 $|$  200 & 100 & 331 &3.81 &  18 &14.09 &  54 &8.51
 \\[3pt]\hline
 qp-rand-m1000-n1500-N100-t1  &1000 $|$ 1000 &1500 $|$ 1500 & 100 & 361 &312.61 &  18 &17175.81 & / & /
 \\[3pt]\hline
 qp-rand-m100-n200-N150-t1  & 100 $|$  100 & 200 $|$  200 & 150 & 341 &6.56 &  19 &20.94 &  58 &60.17
 \\[3pt]\hline
 qp-rand-m1000-n1500-N150-t1  &1000 $|$ 1000 &1500 $|$ 1500 & 150 & 448 &559.20 &  17 &36591.57 & / & /
 \\[3pt]\hline
 qp-rand-m10-n20-N10-t2  &  10 $|$   10 &  20 $|$   20 &  10 &1501 &1.20 &  14 &0.25 & * & *
 \\[3pt]\hline
 qp-rand-m50-n80-N10-t2  &  50 $|$   50 &  80 $|$   80 &  10 & 141 &0.20 &  14 &0.44 & * & *
 \\[3pt]\hline
 qp-rand-m1000-n1500-N10-t2  &1000 $|$ 1000 &1500 $|$ 1500 &  10 & 131 &50.81 &  12 &6916.43 & * & *
 \\[3pt]\hline
 qp-rand-m100-n200-N100-t2  & 100 $|$  100 & 200 $|$  200 & 100 &  81 &3.61 &  14 &28.40 & * & *
 \\[3pt]\hline
 qp-rand-m1000-n1500-N100-t2  &1000 $|$ 1000 &1500 $|$ 1500 & 100 & 220 &576.62 &  13 &8823.43 & * & *
 \\[3pt]\hline
 qp-rand-m100-n200-N150-t2  & 100 $|$  100 & 200 $|$  200 & 150 &  74 &5.36 &  15 &45.91 & * & *
 \\[3pt]\hline
 qp-rand-m1000-n1500-N150-t2  &1000 $|$ 1000 &1500 $|$ 1500 & 150 & 252 &930.81 &  13 &15299.33 & * & *
 \\[3pt]\hline
 qp-tripart1  &2096 $|$  192 &2096 $|$ 2096 &  16 & 653 &1.44 &  15 &1.17 &  24 &0.22
 \\[3pt]\hline
 qp-tripart2  &8432 $|$  768 &8432 $|$ 8432 &  16 & 971 &9.79 &  19 &6.66 &  38 &1.36
 \\[3pt]\hline
 qp-tripart3  &16380 $|$ 1200 &16380 $|$ 16380 &  20 &1034 &27.09 &  22 &30.08 &  55 &6.78
 \\[3pt]\hline
 qp-tripart4  &24815 $|$ 1050 &24815 $|$ 24815 &  35 &5871 &413.35 &  22 &238.92 &  67 &17.46
 \\[3pt]\hline
 qp-gridgen1  &3072 $|$ 1025 &3072 $|$ 3072 & 320 &4081 &308.05 &  40 &2143.19 & 208 &1197.24
 \\[3pt]\hline
 qp-pds15  &1812 $|$ 2125 &7756 $|$ 7756 &  11 &1110 &10.40 &  48 &14.49 &  90 &11.56
 \\[3pt]\hline
 qp-pds30  &3491 $|$ 4223 &16148 $|$ 16148 &  11 &1941 &57.58 &  53 &59.95 & 113 &44.32
 \\[3pt]\hline
 qp-pds60  &6778 $|$ 8423 &33388 $|$ 33388 &  11 &4685 &337.13 &  58 &226.56 & 134 &192.85
 \\[3pt]\hline
 qp-pds90  &8777 $|$ 12186 &46161 $|$ 46161 &  11 &3021 &318.43 &  58 &402.51 & 165 &547.32
 \\[3pt]\hline
 qp-L2CTA3D\textunderscore 100x100x1000\textunderscore 5000  &110000 $|$ 1000 &   0 $|$ 100000 & 100 &  21 &31.24 &   8 &6696.47 &   7 &22.72
 \\[3pt]\hline
 qp-SDC-r100-c50-l100-p1000-t1  &5000 $|$  150 &   0 $|$ 5000 & 100 &  32 &1.52 &   8 &101.91 &   7 &0.84
 \\[3pt]\hline
 qp-SDC-r100-c50-l100-p1000-t2  &5000 $|$  150 &   0 $|$ 5000 & 100 &  31 &1.34 &   6 &96.47 &   6 &0.80
 \\[3pt]\hline
 qp-SDC-r100-c50-l100-p5000-t1  &5000 $|$  150 &   0 $|$ 5000 & 100 &  32 &1.40 &   8 &100.84 &   8 &0.93
 \\[3pt]\hline
 qp-SDC-r100-c50-l100-p5000-t2  &5000 $|$  150 &   0 $|$ 5000 & 100 &  31 &1.37 &   6 &106.82 &   6 &0.79
 \\[3pt]\hline
 qp-SDC-r100-c50-l100-p10000-t1  &5000 $|$  150 &   0 $|$ 5000 & 100 &  32 &1.37 &   8 &97.74 &   8 &0.94
 \\[3pt]\hline
 qp-SDC-r100-c50-l100-p10000-t2  &5000 $|$  150 &   0 $|$ 5000 & 100 &  31 &1.35 &   6 &102.74 &   6 &0.77
 \\[3pt]\hline
 qp-SDC-r100-c100-l100-p1000-t1  &10000 $|$  200 &   0 $|$ 10000 & 100 &  31 &2.67 &   8 &810.58 &   7 &2.16
 \\[3pt]\hline
 qp-SDC-r100-c100-l100-p1000-t2  &10000 $|$  200 &   0 $|$ 10000 & 100 &  31 &2.67 &   6 &1107.95 &   6 &2.10
 \\[3pt]\hline
 qp-SDC-r100-c100-l100-p5000-t1  &10000 $|$  200 &   0 $|$ 10000 & 100 &  31 &2.70 &   8 &1266.75 &   7 &2.12
 \\[3pt]\hline
 qp-SDC-r100-c100-l100-p5000-t2  &10000 $|$  200 &   0 $|$ 10000 & 100 &  31 &2.63 &   6 &751.24 &   6 &2.02
 \\[3pt]\hline
 qp-SDC-r100-c100-l100-p10000-t1  &10000 $|$  200 &   0 $|$ 10000 & 100 &  32 &2.65 &   8 &779.24 &   8 &2.42
 \\[3pt]\hline
 qp-SDC-r100-c100-l100-p10000-t2  &10000 $|$  200 &   0 $|$ 10000 & 100 &  31 &2.63 &   6 &810.03 &   6 &1.98
 \\[3pt]\hline
 qp-SDC-r100-c100-l200-p20000-t1  &10000 $|$  200 &   0 $|$ 10000 & 200 &  41 &6.68 &   8 &1418.31 &   8 &4.87
 \\[3pt]\hline
 qp-SDC-r200-c100-l200-p20000-t1  &20000 $|$  300 &   0 $|$ 20000 & 200 &  34 &11.96 &   8 &5194.45 &   8 &9.47
 \\[3pt]\hline
 qp-SDC-r200-c200-l200-p20000-t1  &40000 $|$  400 &   0 $|$ 40000 & 200 &  31 &22.33 &   8 &53964.31 &   7 &23.34
 \\[3pt]\hline
 qp-SDC-r500-c50-l500-p50000-t1  &25000 $|$  550 &   0 $|$ 25000 & 500 &  41 &43.36 &   8 &11025.98 &   8 &24.56
 \\[3pt]\hline
 qp-SDC-r500-c500-l50-p5000-t1  &250000 $|$ 1000 &   0 $|$ 250000 &  50 &  20 &27.04 &   8 &11360.16 & / & /
 \\[3pt]\hline

	\end{longtable}
\end{center}

Table \ref{table-quad} shows that Gurobi is almost always slowest to solve the
test instances in this case, whereas our sGS-ADMM performs almost as efficiently as BlockIP in solving these quadratic primal block angular problems. It is worth noting that our sGS-ADMM method works very well on the large scale randomly generated problems
compared to BlockIP, because for these instances the
matrices $\cA_i$ and $\cQ_i$ are no longer simple identity matrices for which the BlockIP solver can take special advantage of. Also, BlockIP runs out of memory for three of the huge instances
{\tt qp-rand-m1000-n1500-N100-t1}, {\tt qp-rand-m1000-n1500-N150-t1} and {\tt qp-SDC-r500-c500-l50-p5000-t1}.

It is also observed that BlockIP solver could not solve for the \texttt{qp-rand-xxx-t2} problem because it is not designed to cater for solving problems with nondiagonal quadratic objective cost. For these types of problem, our sGS-ADMM algorithm can substantially outperform Gurobi, sometimes by a factor of more than 10.

%%%%%%%%%%%%%%%%%%%%%%%%%%%%
\subsection{Block angular problems with nonlinear convex objective functions}

In this subsection, we perform numerical experiments on optimization problems having nonlinear convex objective functions and primal block angular constraints.
Nonlinear multicommodity flow problems usually arise in transportation and telecommunication. The two most commonly used nonlinear objective functions are:
\begin{eqnarray*}
	h(t) =
	\begin{cases}
		\sum_{i=1}^m f_{\rm Kr}(t_i; {\rm cap}_i), & \text{known as Kleinrock function;}\\[5pt]
		\sum_{i=1}^m f_{\rm BPR}(t_i; {\rm cap}_i,r_i), & \text{known as BPR (Bureau of Public Roads) function,}
	\end{cases}
\end{eqnarray*}
where

\begin{eqnarray*}
 f_{\rm Kr}(\alp; c) = \left\{ \begin{array}{ll}
  \frac{\alp}{c - \alp} & \mbox{if $0\leq \alp < c$},
  \\[3pt]
  +\infty &\mbox{otherwise},
  \end{array} \right.
  \quad
   f_{\rm BPR}(\alp; c, r) = \left\{ \begin{array}{ll}
  r \alp [1+B(\frac{\alp}{c})^\beta] & \mbox{if $\alp\geq0 $},
  \\[3pt]
  +\infty &\mbox{otherwise}.
 \end{array} \right.
\end{eqnarray*}

The Kleinrock function is normally used to model delay in a telecommunication problem; whereas the BPR function is mainly used to model congestion in a transportation problem. Here
${\rm cap}_i$ is the capacity of arc $i$, $r_i$ is the free flow time of arc $i$, and $\beta,B$ are two positive parameters.

Thus in our problem setting, we have
\begin{eqnarray*}
	&& \theta_0(x_0)=h(x_0), \; \theta_i(x_i)=0, \; \forall\; i = 1,...,N,  \\
	&&\cQ_i = 0, \; c_i = 0, \; \forall\; i = 0,...,N, \; \\
	&&\cA_0 = I, \; \cA_i = -I, \; \forall\; i = 1,...,N,\\
	&&\cD_1 = \cD_2 = \cdots = \cD_N \text{ is a node-arc incidence matrix}, \\
	&&b_0 = 0,\; b_i = d_i \; \forall i = 1,...,N \text{ for some demand } d_i \text{ for each commodity } i, \\
	&& \mathcal{K}_i  = \begin{cases}
		[0,{\rm cap}_i], & \text{for Kleinrock function};\\
		\mathbb{R}_+^{n_i}, & \text{for BPR function} .
	\end{cases}
\end{eqnarray*}

Following \cite{BabonneauVial}, the datasets we used are the \textbf{planar} and \textbf{grid} problems, which could be downloaded from \url{http://www.di.unipi.it/optimize/Data/MMCF.html#Plnr}.

\begin{remark}
	In Step 1c of the sGS-ADMM algorithm, we need to update $s_i^{k+1}$ by
	\begin{eqnarray*}
		s_i^{k+1}
		\;=\;
	\frac{1}{\sig}{\rm Prox}_{\sig \theta_i}\big(\sig(-\cQ_i \bar{w}_i^k + \cD_i^* \bar{y}_i^k + g_i^k)\big) -\big(-\cQ_i \bar{w}_i^k + \cD_i^* \bar{y}_i^k + g_i^k \big)
		 \quad i=0,1,\ldots,N.
	\end{eqnarray*}	
	To compute  the proximal mapping for a given $s$:
	\[{\rm Prox}_{\sigma\theta_i}(s)=\arg\min\big\{g(t):=\sigma\theta_i(t) + \frac{1}{2}\norm{t-s}^2\big\},\]
	we can use   Newton's method to solve the equation $\nabla g(t)=0$. In each sGS-ADMM iteration, we
	warm-start  Newton's method by using the quantity already computed in the
	previous iteration to generate $s_i^k$.
	
	Another point to note is that although $s_i^{k+1}$ is not computed exactly,
	the convergence of the sGS-ADMM algorithm is not affected as long as
	$s_i^{k+1}$ is computed to satisfy the admissible accuracy condition
	required in each iteration of the inexact sGS-ADMM method
	developed in \cite{CST}.
\end{remark}

%%%%%%%%%%%%%%%%%%%%%
\subsubsection{Numerical results}

In this subsection, we compare our sGS-ADMM algorithm against BlockIP and IPOPT. IPOPT is one of the state-of-the-art solvers for solving general nonlinear programs. We use the Kleinrock function as our objective function here.

\begin{center}
	\scriptsize
	\begin{longtable}{| l | c | c | c || c | c || c | c || c | c |}
			\caption{\footnotesize
			Comparison of computational results between sGS-ADMM, BlockIP and IPOPT for nonlinear primal block angular problem. A `-' under the column for BlockIP means that the solver encounters memory issue.}
		\label{table-nonlinear}
		\\
		\hline
		\multicolumn{4}{|c||}{} & \multicolumn{2}{c||}{sGS-ADMM} & \multicolumn{2}{c||}{BlockIP} & \multicolumn{2}{c|}{IPOPT}  \\
		\hline
		\multicolumn{1}{|c|}{Data} & \multicolumn{1}{|c|}{$m_0 | m_i$} & \multicolumn{1}{|c|}{$n_0 | n_i$}
		& \multicolumn{1}{|c||}{$N$}  & \multicolumn{1}{|c|}{Iter} &\multicolumn{1}{|c||}{Time(s)}  & \multicolumn{1}{|c|}{Iter} &\multicolumn{1}{|c||}{Time(s)} &\multicolumn{1}{|c|}{Iter} & \multicolumn{1}{|c|}{Time (s)}\\ \hline
		
		\endhead
		
 grid1  &  80 $|$   24 &  80 $|$   80 &  50 & 591 &0.66 &  28 &0.13 &  76 &3.10 
 \\[3pt]\hline
 grid3  & 360 $|$   99 & 360 $|$  360 &  50 & 381 &0.53 &  41 &1.60 &  86 &21.30 
 \\[3pt]\hline
 grid5  & 840 $|$  224 & 840 $|$  840 & 100 & 581 &1.95 & - & - &  90 &127.50 
 \\[3pt]\hline
 grid8  &2400 $|$  624 &2400 $|$ 2400 & 500 &4171 &261.73 & 215 &4568.28 &  51 &5027.50 
 \\[3pt]\hline
 grid10  &2400 $|$  624 &2400 $|$ 2400 &2000 &3432 &893.98 & 221 &36035.85 &  14 &5340.39 
 \\[3pt]\hline
 planar30  & 150 $|$   29 & 150 $|$  150 &  92 & 431 &0.44 &  93 &1.59 &  90 &7.55 
 \\[3pt]\hline
 planar80  & 440 $|$   79 & 440 $|$  440 & 543 &1875 &20.07 & - & - & 430 &1400.91 
 \\[3pt]\hline
 planar100  & 532 $|$   99 & 532 $|$  532 &1085 &2614 &70.99 & - & - & 117 &1184.46 
 \\[3pt]\hline

	\end{longtable}	
\end{center}

Table \ref{table-nonlinear} shows that IPOPT is almost always the slowest
to solve the test instances but
it is very robust in the sense that it is able to solve all the test instances to the required accuracy.
It is not surprising for it to perform less efficiently
since it is a general solver for nonlinear programs.

On the other hand, we observed that BlockIP runs into memory issue when solving almost half of the instances. This may be due to the fact that BlockIP uses a preconditioned conjugate gradient (PCG) method and Cholesky factorization to solve the linear systems arising in each iteration of the interior-point method. At some point of the iteration, the PCG method did not converge and the algorithm switches to use a Cholesky factorization to solve the linear system, which causes the out-of-memory error.
%We can also observe that BlockIP is no longer efficient in solving these nonlinear primal block angular problems compared to the convex quadratic case.
Even when the PCG method works well, it might still converge in almost $10$ times slower than our algorithm.
%In the other cases, the algorithm simply runs out of memory.

\section{Conclusion}
In conclusion, we have designed efficient methods for solving convex composite quadratic conic programming problems with a primal block angular structure. Numerical experiments show that our algorithm is especially efficient for large instances with convex quadratic objective functions.
As a future project, we plan to
implement our algorithm for solving semidefinite programming problems with primal block angular structures. Also, it would be ideal to utilize a  good parallel computing and
programming platform to implement the algorithm to realize its full potential.

\section{Acknowledgements}
We would like to thank Professor Jordi Castro for sharing with us his
solver BlockIP so that we are able to evaluate the performance of
 our algorithm more comprehensively. We are also grateful
to him for providing us with the test instances he has taken great effort to generate
over the years when developing BlockIP. Thanks also go to the Optimization Group
at the Department of Computer Science of the University of Pisa for
collecting/generating several suites of test data and making them publicly available.

%%%%%%%%%%%%%%%%%%%%%%%%%%%%%%%%%
%% references
%%

\end{document}